\documentclass{article}

\usepackage{amsmath,amssymb,amsthm}
\usepackage{graphicx}
\usepackage{graphics}
\usepackage{epsfig}
\usepackage{hangcaption}
\usepackage{pictexwd}
\usepackage[all]{xy}
\usepackage{amstext}
\usepackage{mathrsfs}
\usepackage{amsfonts}
\usepackage{latexsym}
\usepackage{hyperref}

\theoremstyle{remark}
\newtheorem{defin}{Definition}[section]
\theoremstyle{theorem}
\newtheorem{thm}[defin]{Theorem}
\newtheorem{cor}[defin]{Corollary}
\newtheorem{lem}[defin]{Lemma}
\newtheorem{prop}[defin]{Proposition}
\newtheorem{conj}[defin]{Conjecture}
\newtheorem{cl}[defin]{Claim}
\theoremstyle{remark}
\newtheorem{rem}[defin]{Remark}
\newtheorem{quest}[defin]{Question}
\newtheorem{prob}[defin]{Problem}
\newtheorem{axi}[defin]{Axiom}
\newtheorem{ex}[defin]{Example}
\newtheorem{notice}[defin]{Note}
\newtheorem{conven}[defin]{Convention}

\newcommand{\theorem}[1]{\begin{thm} #1 \end{thm}}
\newcommand{\theoremname}[2]{\begin{thm}[#1] #2 \end{thm}}
\newcommand{\proposition}[1]{\begin{prop} #1 \end{prop}}

\newcommand{\lemma}[1]{\begin{lem} #1 \end{lem}}
\newcommand{\definition}[1]{\begin{defin} #1 \end{defin}}

\newcommand{\remark}[1]{\begin{rem} #1 \end{rem}}

\newcommand{\corollary}[1]{\begin{cor} #1 \end{cor}}

\newcommand{\lemmaname}[2]{\begin{lem}[#1] #2 \end{lem}}
\newcommand{\definitionname}[2]{\begin{defin}[#1] #2 \end{defin}}

\newcommand{\claim}[1]{\begin{cl} #1 \end{cl}}

\newcommand{\ov}[1]{\overline{#1}}
\newcommand{\wh}[1]{\widehat{#1}}

\DeclareMathOperator{\PL}{PL}
\DeclareMathOperator{\Homeo}{Homeo}
\newcommand{\Homeoo}{\Homeo_+}

\DeclareMathOperator{\rot}{Rot}
\newcommand{\ploi}{\PL_+(I)}
\newcommand{\plos}{\PL_+(S^1)}

\newcommand{\homeooI}{\Homeoo(I)}
\newcommand{\homeoS}{\Homeo(S^1)}
\newcommand{\homeooS}{\Homeoo(S^1)}
\newcommand{\ourGps}{\mathcal{A}}
\DeclareMathOperator{\Fix}{Fix}
\DeclareMathOperator{\Supp}{Supp}

\newcommand{\rwr}{\wr_r}
\newcommand{\uwr}{\wr}

\newcommand{\mk}[1]{{#1}}

\newcommand{\ovc}[1]{\overset{\circ}{#1}}

\title{\textsc{Structure Theorems for Groups of
        \\Homeomorphisms of the Circle}}

\author{Collin Bleak \and Martin Kassabov
\thanks{This work was partially supported by NSF grants DMS 0600244, 0635607, 0900932 and a von Neumann fellowship.}
 \and Francesco Matucci\thanks{
This work is part of the third author's PhD thesis at Cornell University.
The third author gratefully acknowledges the Centre de Recerca Matem\`atica (CRM)
and its staff for the support received during the completion of this work.}
}

\begin{document}

\maketitle

\begin{abstract}
In this partly expository paper, we study the set $\ourGps$ of groups
of orientation-preserving homeomorphisms of the circle $S^1$ which do not
admit non-abelian free subgroups.  We use
 classical results about
homeomorphisms of the circle and elementary
 dynamical methods to
derive various new and old results about the groups in $\ourGps$.  Of
the known results, we include some results from a family of results of
Beklaryan and Malyutin, and we also give a new proof of a theorem of
Margulis.  Our primary new results include a detailed classification of the solvable
subgroups of R. Thompson's group $T$.
\end{abstract}

\section{Introduction}
In this paper we explore properties of groups of orientation
preserving homeomorphisms of the circle $S^1$.  In particular, we use
a close analysis of Poincar\'e's rotation number, together with some
elementary dynamical/analytical methods, to prove ``alternative''
theorems in the tradition of the Tits' Alternative.  
Our main result, Theorem \ref{structureThm}, states that any group of orientation
preserving homeomorphisms of the circle is either abelian or is a 
subgroup of a wreath product whose factors can be described in 
considerable detail.
Our methods and
the resulting Theorem \ref{structureThm} give us sufficient
information to derive a short proof of Margulis' Theorem
on the existence of an invariant probability measure on the
circle in \cite{Margulis}, and to
classify the solvable subgroups of the group of orientation-preserving
piecewise-linear homeomorphisms of the circle, and of its subgroup
R. Thompson's group $T$.

Suppose $G$ is a group of orientation preserving homeomorphisms of
$S^1$.  If one replaces the assumption in Theorem \ref{structureThm}
that $G$ has no non-abelian free subgroups with the assumption that
$G$ is a group for which the rotation number map $\rot:G\to\mathbb{R}$
is a homomorphism, then many of the results within the statement of
Theorem \ref{structureThm} can be found in one form or another in the
related works of Beklaryan
\cite{bekla1,bekla2,bekla3,bekla4,bekla5,bekla6}.  However, the
structure of the extension described by Theorem \ref{structureThm} is
new.

A major stepping stone in the established theory of groups of
homeomorphisms of the circle is the following statement (Lemma
\ref{thm:RotHom} below).  {\it For groups of orientation-preserving
  homeomorphisms of the circle which do not admit non-abelian free
  subgroups, the rotation number map is a homomorphism.}  As alluded
in the next paragraph, we believe that the first proof of the
statement comes as a result of combining a theorem of Beklaryan
\cite{bekla3} with Margulis' Theorem in \cite{Margulis}.

Although we arrived at Lemma \ref{thm:RotHom} independently, our
approach to its proof mirrors that of Solodov from his paper
\cite{solodov84}, which also states a version of the lemma as his
Theorem 2.6.  However, in Solodov's proof of his necessary Lemma 2.4,
he uses a construction for an element with non-zero rotation number
that does not actually guarantee that the rotation number is not zero
(see Appendix \ref{sec:appendix}).
Our own technical Lemmas \ref{thm:product-power-conjugate} and
\ref{throw-off} provide sufficient control to create such an element,
and the rest of the approach goes through unhindered.

As also shown by Beklaryan, the results within Theorem
\ref{structureThm} can be employed to prove Margulis' Theorem.  Our
proof of Theorem \ref{structureThm} (including Lemmas
\ref{thm:product-power-conjugate} and \ref{throw-off}) and the
consequential proof of Margulis' Theorem, both use only classical
methods.

Due to the large intersection with known work and results, portions of
this paper should be considered as expository.  Many of
the proofs we give are new, taking advantage of our technical Lemma
\ref{throw-off}. This lemma may have other applications as well.  
Further portions of this project, which trace out some new proofs of
other well-known results, are given in the third author's
dissertation~\cite{MatucciThesis}.

We would like to draw the reader's attention to the surveys by
Ghys~\cite{GhysSurvey} and Beklaryan~\cite{bekla6} on groups of
homeomorphisms of the circle, and to the book by
Navas~\cite{NavasBook} on groups of diffeomorphisms of the circle, as
three guiding works which can lead the reader further into the
theory.

\subsection*{Statement and discussion of the main results} 
We use much of the remainder of the introduction to state and briefly
discuss our primary results. Except for Lemma
~\ref{thm:RotHom}, Theorem ~\ref{MargulisThm} and parts of Theorem
~\ref{structureThm}, our results are new.

\mk{\subsection{The main structure theorem}}

Denote by $\homeooS$ the maximal subgroup of $\homeoS$ consisting of
orientation preserving homeomorphisms of $S^1$ and let $\rot: \homeooS
\to \mathbb{R}/\mathbb{Z}$ denote Poincar\'e's rotation number
function.  Although this function is not a homomorphism, we will
denote by $\ker \left(\rot\right)$ its ``kernel'', i.e., the set of
elements with rotation number equal to zero.  Similarly, denote by
$\homeooI$ the maximal group of orientation-preserving homeomorphisms
of the unit interval.  

In order to state our first result, we note that by
Lemma~\ref{thm:RotHom} the restriction of Rot to any subgroup of
$\homeooS$ which has no non-abelian free subgroups turns Rot into a
homomorphism of groups.  Also, note that throughout this article we
use the expressions $C\uwr T \simeq \left(\prod_{t\in
  T}C\right)\rtimes T$ and $C\rwr T\simeq \left(\bigoplus_{t\in
  T}C\right)\rtimes T$ respectively to denote the unrestricted and
restricted standard wreath products of groups $C$ and $T$.

\medskip
 \theorem{\label{structureThm}

Let $G\in \homeooS$, with no non-abelian
 free subgroups.  Either $G$
is abelian or there are subgroups $H_0$ and $Q$ of
$\homeooS$, such that
\[
G\hookrightarrow H_0 \uwr Q
\]
 where the embedding is such that the following hold.
\begin{enumerate}
\item The group $H_0$ has the following properties.
\begin{enumerate}
\item $\rot$ is trivial over $H_0$.
\item There is an embedding $H_0\hookrightarrow\prod_{\mathfrak{N}} \homeooI$,
  where $\mathfrak{N}$ is an index set which is at most countable.
\item The group $H_0$ has no non-abelian free subgroups.
\end{enumerate}
\item The group $Q \cong G/\left(\ker\left(\rot\right)\cap G\right)$
  is isomorphic to a subgroup of $\mathbb{R}/\mathbb{Z}$, which is
  at most countable.
\item The subgroups $H_0,Q \le \homeooS$ generate a subgroup
  isomorphic to the restricted wreath product $H_0 \rwr Q$. This
  subgroup can be ``extended'' to an embedding of the unrestricted
  wreath product into $\homeooS$ where the embedded extension contains
  $G$.
\end{enumerate}
}

\remark{We note that if the kernel of the homomorphism $\rot$ is
  trivial over $G$ then $G$ embeds in a pure group of rotations and so
  is abelian.}

As mentioned in the introduction, most of points one and two above can
be extracted from the results of Beklaryan in \cite{bekla2,bekla3,bekla4} under the assumption 
of the existence of a $G$-invariant probability measure on the circle (a property which 
Beklaryan shows to be equivalent to $\rot:G\to \mathbb{R}/\mathbb{Z}$ being a homomorphism in \cite{bekla3}).

Theorem 1.1 attempts to provide an algebraic description of a
dynamical picture painted by Ghys in~\cite{GhysSurvey}.  We will quote
a relevant statement below to clarify this comment.  First though, we
give a description of these same dynamics using the construction of a
counter-example to Denjoy's Theorem in the $C^1$ category (there is a
detailed, highly concrete construction of this counter-example
in~\cite{TopologyFoliations}, and a detailed discussion of a family of
counter-examples along these same lines in section 4.1.4 of
~\cite{NavasBook}).

Denjoy's Theorem states
that given a $C^2$ orientation-preserving
 circle homeomorphism
$f:S^1\to S^1$ with irrational rotation number
 $\alpha$ (in some
sense, points are moved ``on average'' the distance
 $\alpha$ around
the circle by $f$), then there is a homeomorphism
 $c:S^1\to S^1$ so
that $c\circ f\circ c^{-1}$ is a pure rotation of
 the circle by
$\alpha$.

We now discuss the counter-example: Take a rotation $r$ of the circle
by an irrational $\alpha$ ($r$ is a circle map with real lift map
$t\mapsto t+\alpha$, under the projection map $p(t) = e^{2\pi i t}$).
The orbit of any point under iteration of this map is dense on the
circle.  Now, track the orbit of a particular point in the
circle.  For each point in the orbit, replace the point by an interval
with decreasing size (as our index grows in absolute value), so that
the resulting space is still homeomorphic to $S^1$.  Now, extend $r$'s
action over this new circle so that it becomes a $C^1$ diffeomorphism
$\tilde{r}$ of the circle which agrees with the original map $r$ over
points in the original circle, and which is nearly affine while
mapping the intervals to each other.\footnote{This can be done if the
  interval lengths are chosen carefully. However, by Denjoy's Theorem,
  no matter how one chooses lengths and the extension $\tilde{r}$, the
  result will fail to be a $C^2$ diffeomorphism of the circle.}  The
map $\tilde{r}$ still has the same rotation number as $r$, and cannot
be topologically conjugated to a pure rotation because there are points whose
orbits are not dense.
 
 Let $H_0$ be any group of orientation-preserving
 homeomorphisms of
 the interval.  Pick an element of $H_0$ to act on
 one of the
 ``inserted''
 intervals above, and further elements in
 copies of
 $H_0$ (created by
 conjugating the original action of
 $H_0$ by powers
 of $\tilde{r}$) to
 act on the other ``inserted''
 intervals.  We have just constructed an element of $H_0\uwr
 \mathbb{Z}$, acting on (a scaled up version) of $S^1$.

While providing a useful picture, the above explanation does not
really capture the full dynamical picture implied by 
Theorem~\ref{structureThm}; the group $G$ may be any subgroup of the
appropriate wreath product, so elements of the top group in the wreath
product may
 not be available in $G$.  Further, based on possible
categorical restrictions on the group $G$, 
other restrictions on the wreath
product may come into play.

Now let us relate this picture to Ghys' discussion
in~\cite{GhysSurvey}.  
 In a sentence near the end of the final
paragraph
 of section 5 in~\cite{GhysSurvey}, Ghys states the
following.
 
\begin{quote}
{\it $\ldots$ we deduce that [$G$] contains a non abelian free
  subgroup unless the restriction of the action of [$G$] to the
  exceptional minimal set is abelian and is semi-conjugate to a group
  of rotations $\ldots$ }
\end{quote}
Here, the complement of the exceptional minimal set of the action
of $G$ contains the region where our base group acts, and the
top group acts essentially as (is semi-conjugate to) a group of
rotations on the resultant circle which arises after ``gluing
together'' the exceptional minimal set (using the induced cyclic
ordering from the original circle).

\subsection{Some embedding theorems}
The theorems in this subsection follow by combining the results
(see~\cite{bpasc,bpnsc, navas1}) of the first author or of Navas on
groups of piecewise-linear homeomorphisms of the unit interval
together with Theorem \ref{structureThm}.

Throughout this article, we will use $\ploi$ and $\plos$ to represent
the piecewise-linear orientation-preserving homeomorphisms of the unit
interval $I:=[0,1]$ and of the circle $S^1$, respectively.

In order to state our embedding results and to trace them as
consequences of Theorem~\ref{structureThm}, we need to give some
definitions and results from~\cite{bpasc,bpnsc}.  Let $G_0=1$ and, for
$n \in \mathbb{N}$, inductively define $G_n$ as the direct sum of a
countably infinite collection of copies of the group $G_{n-1} \rwr
\mathbb{Z}$:
\[
G_n := \bigoplus_{\mathbb{Z}} \left(G_{n-1} \rwr \mathbb{Z} \right).
\]

A result in~\cite{bpasc} states that if $H$ is a solvable group
with derived length $n$, then $H$ embeds in $\ploi$ if and only
if $H$ embeds in $G_n$.  Using Theorem~\ref{structureThm} and 
Remark~\ref{thm:PL-is-always-rational} (see section 5), we are able to extend
this result to subgroups of $\plos$:

\medskip
\theorem{
\label{thm:solvClassPL}
Suppose $H$ is a solvable group with derived length $n$.  The group $H$ embeds in $\plos$ if and only if one of the following holds,
\begin{enumerate}
\item $H$ embeds in $\mathbb{R}/\mathbb{Z}$,
\item $H$ embeds in $G_n$, or
\item $H$ embeds in $G_{n-1} \rwr K$ for some nontrivial subgroup $K$ of
  $\mathbb{Q}/\mathbb{Z}$.

\end{enumerate}
}

The paper~\cite{bpnsc} also gives a non-solvability criterion for
subgroups of $\ploi$. Let $W_0=1$ and, for $n \in \mathbb{N}$,
we define $W_i= W_{i-1} \rwr \mathbb{Z}$.  Build the group
\[
W := \bigoplus_{i \in \mathbb{N}} W_i.
\]

The \mk{main} 
result of~\cite{bpnsc} is that a subgroup $H\le \ploi$ is
non-solvable if and only if $W$ embeds in $H$.  Now again by using
Theorem~\ref{structureThm}, we are able to give a Tits' Alternative
type of theorem for subgroups of $\plos$:

\theorem
{\label{thm:nonSolvClassPL}
A subgroup $H \le\plos$ either
\begin{enumerate}
\item contains a non-abelian free subgroup on two generators, or
\item contains a copy of $W$, or
\item is solvable.
\end{enumerate}
}

As may be clear from the discussion of the counterexample to Denjoy's Theorem, 
it is not hard to produce various required wreath products as groups of homeomorphisms of the circle.
\theorem{
\label{thm:first-embedding}
For every countable subgroup $K$ of $\mathbb{R} / \mathbb{Z}$ and for
every $H_0 \le \homeooI$ there is an embedding $H_0 \uwr K
\hookrightarrow \homeooS$. }

We recall the R. Thompson groups
$F$ and $T$.  These are groups of
 homeomorphisms of the interval $I$
and of the circle
 $\mathbb{R}/\mathbb{Z}$ respectively.  In
particular, they are the groups one obtains if one restricts the
groups of orientation preserving homeomorphisms of these spaces to the
piecewise-linear category, and insist that these piecewise linear
elements (1) have all slopes as integral powers of two, (2) have all
changes in slope occur at dyadic rationals, and (3) map the dyadic
rationals to themselves. 

\theorem{
\label{thm:embedding-thompson}
For every $K \le \mathbb{Q} / \mathbb{Z}$ there is an embedding $F
\rwr K \hookrightarrow T$, where $F$ and $T$ are the
R. Thompson groups above.}

More generally, we have the following similar theorem.

\theorem{
\label{thm:embedding-PL}
For every $K \le \mathbb{Q} / \mathbb{Z}$
there is an embedding $\ploi \rwr K \hookrightarrow \plos$. }

\subsection{Useful Lemmas}
 Our proof of the following lemma sets the foundation upon which the
 other results in this article are built.  As mentioned in the
 introduction, the standing proof of Lemma~\ref{thm:RotHom} is to
 quote Theorem 6.7 of~\cite{bekla3}, together with Margulis' Theorem
 (Theorem~\ref{MargulisThm} below).  \lemmaname{Beklaryan and Margulis,~\cite{bekla3, Margulis}} {
\label{thm:RotHom}
Let $G \le \homeooS$. Then the following alternative holds:
\begin{enumerate}
\item $G$ has a non-abelian free subgroup, or
\item the map $\rot:G \to (\mathbb{R}/\mathbb{Z},+)$ is a group homomorphism.
\end{enumerate}
}
 
The heart of the proof of Lemma~\ref{thm:RotHom} is contained in the
following lemma, which itself is proven using only on classical
results (Poincare's Lemma and the Ping-pong Lemma).  We mention the
lemma below in this section as it provides a useful new technical
tool.

In the statement below, if $G$ is a group of homeomorphisms of the
circle, and $g\in G$, then Fix($g$) is the set of points of the
circle
 which are fixed by the action of $g$ and $G_0=\left\{g\in G
\mid
 \Fix(g) \ne \emptyset\right\}$.

\lemmaname{Finite Intersection Property}
{\label{thm:FIP}
Let $G \le \homeooS$ with no non-abelian free subgroups.
The family $\{\Fix(g) \mid  \, g \in G_0\}$ 
satisfies the finite intersection property,
i.e., for all $n$-tuples $g_1, \ldots,g_n \in G_0$, we have
$\Fix(g_1) \cap \ldots \cap \Fix(g_n) \ne \emptyset$.}

Another view of the above lemma is \mk{the following ``generalization'' of the
Ping-pong lemma: let $X$ be a collection}
of homeomorphisms of the circle such that
\begin{enumerate}
\item for all $g\in X$, Fix($g$) $\neq \emptyset$, and
\item for all $x\in S^1$ there is some $g\in X$ with $g(x) \neq x$,
\end{enumerate}
then $\langle X\rangle$ contains embedded non-abelian free groups.

\mk{\subsection{Some further applications}}
As mentioned above, our proof of Lemma~\ref{thm:RotHom} uses only
elementary methods and classical results.  Margulis' Theorem follows
very simply with Lemma~\ref{thm:RotHom} in hand.  We
hope our approach
provides a valuable new perspective on this theorem.

\medskip
\theoremname{Margulis,~\cite{Margulis}}{
\label{MargulisThm}
Let $G \le \homeooS$. Then at least one of the two following statements must be
true:
\begin{enumerate}
\item $G$ has a non-abelian free subgroup, or
\item there is a $G$-invariant probability measure on $S^1$.
\end{enumerate}
}

Finally, we mention a theorem which gives an example of how
restricting the category gives added control on the wreath product of
the main structure theorem.  It may be that the following result is
known, but we were not able to find a reference for it.  The following
application represents the only occasion where we rely upon Denjoy's Theorem.

\theorem{
\label{thm:strongTits}
\label{thm:irrational-abelian}
Suppose $G$ is a subgroup of $\homeooS$ so that
the elements of $G$ are either
\begin{enumerate}
\item all piecewise-linear, each admitting at most finitely many
  breakpoints, or
\item all $C^1$ with bounded variation in the first derivative,
\end{enumerate}
\item 
 and suppose there is $g\in G$ with $\rot(g) \not \in
  \mathbb{Q} / \mathbb{Z}$.
 Then $G$ is topologically conjugate to a
  group of rotations (and is thus abelian) or $G$ contains a
  non-abelian free subgroup.}

\subsection*{Organization}
The paper is organized as follows: Section~\ref{sec:back-tools}
recalls the necessary language and tools which will be used in the
paper; Section~\ref{sec:rot-homo} shows that the rotation number map is a
homomorphism under certain hypotheses; Section~\ref{sec:margulis}
uses the fact that the rotation number map is a homomorphism to prove
Margulis' Theorem on invariant measures on the unit circle; 
Section~\ref{sec:structure-embedding} proves and demonstrates the main
structure theorem.  
\subsection*{Acknowledgments}
The authors would like to thank Ken Brown and Matt Brin for a careful reading
of some versions of this paper and many helpful remarks. The
authors would also like to thank Isabelle Liousse for helpful
conversations and L.A. Beklaryan, \'Etienne Ghys, Andr\'es Navas and Dave
Witte-Morris for providing important references. 
We would also like to thank various referees and editors whose comments have helped to develop this article and suggested that we add an appendix to clarify a statement.

We owe a particular debt of gratitude to Serge Cantat, a visitor
 at
Cornell University at the time of this research, who informed
 us of
Margulis' Theorem, and suggested that we try to find a proof of it
from our own point of view.
 
\section{Background and Tools}
\label{sec:back-tools}

 In this section we collect some known results we will use throughout
the paper.  We use the symbol $S^1$ to either represent
$\mathbb{R}/\mathbb{Z}$ (in order to have a well defined origin $0$)
or as the set of points in the complex plane with distance one from
the origin, as is convenient.  We begin by recalling the definition
of rotation number. Given $f \in \homeooS$, let $F: \mathbb{R} \to
\mathbb{R}$
 represent a lift of $f$ via the standard covering
projection $\exp : \mathbb{R} \to S^1$,
 defined as $\exp(t) = e^{2
\pi it}$.

Following~ \cite{poincare1, poincare2}, we  define the rotation number of an
orientation-preserving homeomorphism of the circle. 
Consider the limit
\begin{equation}
\label{eq:rotation}
\lim_{n \to \infty} \frac{F^n(t)}{n} \pmod{1}.
\end{equation}
It is possible to prove that this limit exists and that it is
independent of the choice of $t$ used in the above calculation 
(see~\cite{Hermann}).
 Moreover, such a limit is independent of the choice
of lift \mk{$F$, when considered} $\pmod{1}$.

\definitionname{Rotation number of a function}{Given $f \in \homeooS$
and $F \in \Homeo(\mathbb{R})$ a lift of $f$, we say that
\[
\lim_{n \to \infty} \frac{F^n(t)}{n} \pmod{1} := \rot(f) \in  \mathbb{R}/\mathbb{Z}
\]
is the \emph{rotation number} of $f$.}

\definition{Given $f \in \homeooS$, we define $\Fix(f)$ to be the set of points
that are fixed by $f$, i.e. $\Fix(f)=\{s \in S^1 \mid f(s)=s\}$. A similar definition
is implied for any $F \in \Homeoo(\mathbb{R})$.}

Since the rotation number is independent of the choice of the lift, we will work with
a preferred lift of elements and of functions.

\definitionname{The ``hat'' lift of a point and of a function}{
\label{thm:hat-lift-def}
For any element $x \in S^1$ we denote by $\wh{x}$ the lift of $x$ contained in $[0,1)$.
For functions in $\homeooS$ we distinguish between functions with or without fixed points
and we choose a lift that is ``closest'' to the identity map. If
$f \in \homeooS$ and the fixed point set $\Fix(f) = \emptyset$, we denote by $\wh{f}$
the lift to $\Homeoo(\mathbb{R})$ such that $t < \wh{f}(t) < t+1$ for all $t \in \mathbb{R}$.
If $f \in \homeooS$ and $\Fix(f) \ne \emptyset$,
we denote by $\wh{f}$ the lift to $\Homeoo(\mathbb{R})$ such that $\Fix(\wh{f}) \ne \emptyset$. 
The map $\wh{f}$ can also be defined as the unique lift such
that $0\leq \lim_{n \to \infty} \frac{\wh{f}^n(t)}{n} < 1$, for all $t \in \mathbb{R}$.
}

We will use these definitions for lifts of elements and functions in Lemma~\ref{thm:lift-control}(4) and
throughout the proof of Lemma~\ref{thm:RotHom}.
If we use this lift to compute the limit defined in (\ref{eq:rotation}), the result is always in $[0,1)$.
Proofs of the next three results can be found in~\cite{Hermann} and~\cite{mackay}.

\lemmaname{Properties of the Rotation Number}{
\label{thm:lift-control}
Let $f,g \in \homeooS$, $G\le \homeooS$ and $n$ be a positive
integer. Then:
\begin{enumerate}
\item $\rot(f^g)=\rot(f)$.

\item $\rot(f^n)=n \cdot \rot(f)$.

\item If $G$ is abelian then the map
\[
\begin{array}{cccc}
\rot: & G & \longrightarrow & \mathbb{R}/\mathbb{Z} \\
     & f & \longmapsto     & \rot(f)
\end{array}
\]
is a group homomorphism.

\item If $\rot(g)=p/q \pmod{1} \in \mathbb{Q}/\mathbb{Z}$ and $s \in S^1$ is such that $g^q(s)=s$,
then $\wh{g}^q(\wh{s})=\wh{s}+p$.
\end{enumerate}
}

Two of the most important results about the rotation number are stated below:

\theoremname{Poincar\'e's Lemma}{
Let $f \in \homeooS$ be a homeomorphism. Then
\begin{enumerate}
\item $f$ has a periodic orbit of length $q$ if and only if $\rot(f)=p/q \pmod{1} \in \mathbb{Q}/\mathbb{Z}$ and
$p,q$ are coprime.
\item $f$ has a fixed point if and only if $\rot(f) =0$.
\end{enumerate}}

\medskip

We recall that \emph{Thompson's group $T$} is the subgroup of
elements of $\plos$ such that for any such element all breakpoints
occur at dyadic rational points, all slopes are powers of $2$, and
dyadic rationals are mapped to themselves.  Moreover, recall that the
subgroup of $T$ consisting of all elements which fix the origin $0$ is
one of the standard representations of \emph{Thompson's group $F$}
(for an oft-cited introduction about Thompson's groups, see~\cite{cfp}). 
Ghys and Sergiescu prove in~\cite{GhysSergiescu} that
all the elements of Thompson's group $T$ have rational rotation
number. Liousse in~\cite{Liousse1} generalizes this result to the
family of \emph{Thompson-Stein groups} which are subgroups of $\plos$
with certain suitable restrictions on rational breakpoints and slopes.

\medskip

The following is a classical result proved by Fricke and Klein
\cite{FK} which we will need in the proofs of section~\ref{sec:rot-homo}.


\theoremname{Ping-pong Lemma}{
\label{thm:ping-pong}
Let $G$ be a group of permutations on a
  set $X$, let $g_1, g_2$
 be elements of $G$.
If $X_1$ and $X_2$ are disjoint subsets of $X$ and for all
  integers
 $n \ne 0$, $i \ne j$, $g_i^n(X_j) \subseteq X_i$ , then
  $g_1, g_2$ freely generate the free group $F_2$
 on two
  generators. }

We use the following theorem only to give an application of our main
structure theorem.  The version we give below is an expansion of
Denjoy's original theorem.  An elegant proof of the content of this
statement is contained in the paper~\cite{mackay}.

\theoremname{Denjoy~\cite{Denjoy}}{
\label{thm:denjoy}
Suppose $f \in \homeooS$ is
piecewise-linear with finitely many breakpoints
or is a $C^1$
homeomorphism whose first derivative has bounded variation. If the
rotation number of $f$ is irrational, then $f$ is conjugate (by an
element in $\homeooS$) to a rotation. Moreover, every orbit of
$f$ is dense in $S^1$. }

\section{The Rotation Number Map is a Homomorphism}
\label{sec:rot-homo}
Our main goal for this section is to prove 
Lemma~\ref{thm:RotHom}, which states that the rotation number map is a homomorphism under certain assumptions.
It is not true in general that the rotation number map is a group homomorphism. 
The example drawn in figure
\ref{fig:rot-number-zero-not-homomorphism-4} below shows a pair of maps with fixed points (hence with
rotation number equal to zero, by Poincar\'e's Lemma) and such that their product does not fix any point
(thus has non-zero rotation number).

\begin{figure}[0.5\textwidth]
 \begin{center}
  \includegraphics[height=4.5cm]{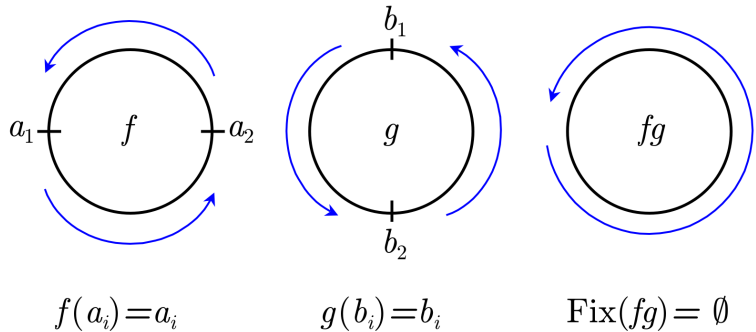}
 \end{center}
 \caption{The rotation number map is not a homomorphism in general.}
 \label{fig:rot-number-zero-not-homomorphism-4}
\end{figure}

\definition{
We define the (open) \emph{support} of $f$ to be the set of points which are moved by $f$, i.e.,
$\Supp(f)=S^1 \setminus \Fix(f)$.\footnote{Notice that this definition is a bit different
from the definition in analysis, where supports are forced to be closed sets.} 
 A similar definition is implied for any $f \in \Homeoo(\mathbb{R})$. }

Our proof divides naturally into several steps. We start by showing
how to use the Ping-pong Lemma
 to create free subgroups. This idea
is well known (see for example Lemma 4.3 in~\cite{brin1}), but we give
an account of it for completeness.

\lemma{
\label{thm:free-group}
Let $f,g \in \homeooS$ such that $\Fix(f) \ne \emptyset \ne \Fix(g)$.
If the intersection $\Fix(f) \cap \Fix(g) = \emptyset$, 
then $\langle f,g \rangle$ contains a
non-abelian free subgroup. 
}

\begin{proof}
Let $S^1 \setminus \Fix(f)=\bigcup I_\alpha$ and $S^1 \setminus \Fix(g)=\bigcup J_\beta$,
for suitable families of pairwise disjoint open intervals $\{I_\alpha\},\{J_\beta\}$.
We assume $\Fix(f) \cap \Fix(g) = \emptyset$
 so
that $S^1 \subseteq \left(\bigcup I_\alpha\right) \cup \left(\bigcup
J_\beta\right)$.

Since $S^1$ is compact, we can write
$S^1=I_1\cup \ldots \cup I_r \cup J_1 \cup \ldots \cup J_s$, for suitable
intervals in the families $\{I_\alpha\},\{J_\beta\}$.
 Define $I=I_1
\cup \ldots \cup I_r$ and $J=J_1 \cup \ldots \cup J_s$.
We observe that $\partial I$ and $\partial J$ are finite and that,
since each $x \in \partial J$ lies in the interior of $I$, there is an open
neighborhood $U_x$ of $x$ such that $U_x \subseteq I$.
 Let $X_g =
\bigcup_{x \in \partial J} U_x$. Similarly we build an open set
$X_f$. The neighborhoods used
 to build $X_f$ and $X_g$ can be chosen
to be small enough so that $X_f \cap X_g =\emptyset$.  If $x \in
\partial J$,
 then the sequence $\{f^n(x)\}_{n \in \mathbb{N}}$
accumulates at a point of $\partial I$ and so there is
 an $n \in
\mathbb{N}$ such that $f^n(U_x) \subseteq X_f$.
 By repeating this
process for each $x \in \partial J$ and $y \in \partial I$,
we find an $N$ big enough so that for all $m \ge N$ we have
\[
f^m(X_g) \cup f^{-m}(X_g) \subseteq X_f, \; \; g^m(X_f) \cup g^{-m}(X_f) \subseteq X_g.
\]
If we define $g_1 = f^N,g_2=g^N,X_1=X_f,X_2=X_g$, we satisfy the
hypothesis of Theorem~\ref{thm:ping-pong} since both of the elements
$g_1,g_2$ have infinite order.
 Thus $\langle g_1,g_2 \rangle$ is a
non-abelian free subgroup of $\langle f,g\rangle$.
\end{proof}

\corollary{
\label{thm:free-group-2}
Let $f,g \in \homeooS$ such that $\Fix(\wh{f}) \ne \emptyset \ne \Fix(\wh{g})$.
If $\Fix(\wh{f}) \cap \Fix(\wh{g}) = \emptyset$, 
then $\langle f,g \rangle$ contains a non-abelian free subgroup.
}

\definition{
If $G \le \homeooS$ is a group, as in the introduction we define the set of homeomorphisms with fixed points
\[
G_0=\{g \in G \mid \exists s\in S^1, g(s)=s\} = \{ g \in G \mid \rot(g)=0 \} \subseteq G.
\]}

\corollary{
Let $G \le \homeooS$ with no non-abelian free subgroups.
The subset $G_0$ is a normal subgroup of $\homeooS$.
}

\begin{proof}
Let $f,g \in G_0$ then, by Lemma~\ref{thm:free-group}, they must have a common
fixed point, hence $fg^{-1} \in G_0$ and $G_0$ is a subgroup of $G$.
Moreover, if $f \in G,g \in G_0$
and $s \in \Fix(g)$, we have that $f^{-1}(s) \in \Fix(f^{-1}gf)$ and so that $f^{-1}gf \in G_0$
and therefore $G_0$ is normal.
\end{proof}

If $f$ has no fixed points then the support of $f$ is the whole
circle $S^1$, otherwise the support can be broken into%
\mk{\footnote{possibly infinitely many}} open
intervals upon each of which $f$ acts as a one-bump function, 
that is $f(x) \ne x$ on each such interval.

\definition{
Given $f \in \homeooS$, we define an \emph{orbital of $f$} 
as a connected component of the support of $f$.  
If $G\leq \homeooS$ then we define an \emph{orbital of $G$} as a connected
component of the support of the action of $G$ on $S^1$.
}

We note in passing that any orbital of $G$ can be written as a union of orbitals
of elements of $G$.

Lemmas \ref{thm:build-large}, \ref{thm:product-power-conjugate} , 
and \ref{throw-off} are highly technical lemmas from which one easily
derives the useful Corollary \ref{thm:compactness-global}. While Lemmas \ref{thm:build-large}--\ref{throw-off} are proven using
elementary techniques, these Lemmas and the techniques involved in their
proofs have no bearing on the remainder of the paper.  Thus, the reader more
interested in the global argument will not lose much by passing directly to
Corollary \ref{thm:compactness-global} on an initial reading.

The following lemma is straightforward and can be derived using
techniques similar to those of the first author in~\cite{bpgsc} or
those of Brin and Squier in~\cite{brin1}. We omit its proof.

\lemma{
\label{thm:build-large}
Let $H\le \homeooI$ and let $(a,b)$ be an interval
such that $\Fix(H) \cap (a,b)=\emptyset$.
For every $\varepsilon >0$, there is an element
$w \in H$ such that $w$ has an orbital
containing $[a+\varepsilon, b-\varepsilon]$.
}

The following will be used in the proof of 
Lemma~\ref{throw-off}.

\lemma{
\label{thm:product-power-conjugate} 
Let $H \le \homeooI$ and suppose that $(a_1,b_1),\ldots,(a_r,b_r)$
are orbitals of $H$.  Let $\varepsilon>0$ and suppose 
there is an element $f \in H$ such that 
$\Supp(f) \supseteq \cup [a_i+\varepsilon,b_i-\varepsilon]$.
Given any $g \in H$ there exists a positive integer $M$ such that for all $m \ge M$, 
there exist positive integers $K$ 
and $N$ 
such that for all $n \ge N$, we have
\[
f^m g^n f^{-m} \cdot f^{-K}\, \left(\bigcup_{i=1}^{r}[a_i+\varepsilon,b_i-\varepsilon]\right)\,\cap\, \bigcup_{i=1}^{r}[a_i+\varepsilon,b_i-\varepsilon] \,=\, \emptyset.
\]
}

\begin{proof}
We consider the set $\mathcal{J}=\{(s_1,t_1),\ldots,(s_r,t_r)\}$ of components
of the support of $f$ respectively containing the intervals 
$[a_i+\varepsilon,b_i-\varepsilon]$.

Fix an index $i$ and let us suppose for now that $f(x)>x$ for all
$x\in (s_i,t_i)$.
We consider the possible fashions in which $g$ can have support in 
$(a_i,b_i)$, where the actions of $g$ and $f$ may interact.

There are three cases of interest.

\begin{enumerate}
\item There is an orbital $(u_i,v_i)$ of $g$ such that $s_i\in [u_i,v_i)$.
\item There is a non-empty interval $(s_i,x_i)$ upon which $g$ acts as the identity.
\item The point $s_i$ is an accumulation point of a decreasing sequence 
      of left endpoints $\{x_{i,j}\}_{j\in \mathbb{N}}$ of orbitals of $g$ contained in $(s_i,t_i)$.
\end{enumerate}

In the first case, since $f$ is increasing on $(s_i,t_i)$, there
exists a positive power $M_i$ such that $f^m(v_i)>b_i-\varepsilon$ for
all $m\ge M_i$.  Hence any such conjugate $f^m g f^{-m}$ will have an
orbital containing $(s_i,b_i-\varepsilon]$.  For any $K_i>0$ the set
  $W_{i,K_i}:=f^{-K_i}([a_i+\varepsilon,b_i-\varepsilon])$ is a
  compact connected set inside $(s_i,b_i-\varepsilon)$, hence there
  exists an $N_i>0$ such that for all $n>N_i$ we have
\[
f^m g^n f^{-m}(W_{i,K_i}) \cap [a_i+\varepsilon,b_i-\varepsilon] = \emptyset.
\]

In the second case we assume that $g$ is the identity on an interval
$(s_i,x_i)$, for some $s_i<x_i<b_i$.  There exists a non-negative
power $M_i$ such that $f^m(x_i)>b_i-\varepsilon$ for all $m\ge M_i$.
Hence any conjugate $f^m g f^{-m}$ for $m\ge M_i$ will be the identity
on the interval $(s_i,b_i-\varepsilon]$.  In particular, if $K$ is
large enough so that $W_{i,K}\cap [a_i +\varepsilon, b_i -
  \varepsilon] = \emptyset$ we must have that for any $m > M_i$ and
any integer $n \ge N_i$ (for any positive integer $N_i$) the product
$f^m g^n f^{-m}f^{-K_i}$ will move $[a_i +\varepsilon, b_i -
  \varepsilon]$ entirely off of itself.

In the third case we assume that $s_i$ is the accumulation point of a
decreasing sequence of left endpoints $\{x_{i,j}\}_{j \in \mathbb{N}}$
of orbitals of $g$ contained in $(s_i,t_i)$.  Given any positive
integer $M_i$ observe that if $m \ge M_i$, there exists an index $j_m$
such that $x_m:=f^m(x_{i,j_m})<a_i+\varepsilon$.  Let $N_i = 1$ and
note that for any power $n\geq N_i$ the conjugate $f^m g^n f^{-m}$ fixes
$x_m$. Now we choose $K_i$ to be large enough so that
$f^{-K_i}(b_i-\varepsilon)<x_m$. With these choices, the product $f^m
g^n f^{-m} \cdot f^{-K_i}$ moves $[a_i+\varepsilon,b_i-\varepsilon]$
entirely off of itself.

We note in passing that in all three cases, $K_i$ could always be
chosen larger, with the effect (and only in the first case) that we
might have to choose $N_i$ larger.

If instead $f$ is decreasing on the interval $(s_i,t_i)$, similar
(reflecting right and left) arguments based at the point $t_i$ instead
of $s_i$, will find products $f^m g^n f^{-m} \cdot f^{-K_i}$ in all
three corresponding cases which move
$[a_i+\varepsilon,b_i-\varepsilon]$ entirely off of itself.

Choose $M=\max\{M_1,\ldots,M_r\}$ and choose any $m \ge M$.  Given
this choice of $m$ there are minimal positive choices of $K_i$ for
each index $i$ as above.  Set $K=\max\{K_1,\ldots,K_r\}$.  For this
choice of $K$ we can find, for each index $i$, an integer $N_i$ so
that for all values of $n> N_i$, our product will move
$[a_i+\varepsilon,b_i -\varepsilon]$ entirely off of itself.  Now set
$N=\max\{N_1,\ldots,N_r\}$. With these choices, we have that for all
$n\ge N$ the product $f^m g^n f^{-m} \cdot f^{-K}$ moves every
$[a_i+\varepsilon,b_i-\varepsilon]$ entirely off of itself for all
indices $i$.
\end{proof}

\lemma{
\label{throw-off}
Let $H \le \homeooI$, let $(a_1,b_1),\ldots,(a_r,b_r)$
be a finite collection of components of the support of $H$, and let $\varepsilon>0$. Then
there exists $w_\varepsilon \in H$ such that for all $i$
\begin{equation}
w_\varepsilon\left([a_i+\varepsilon,b_i-\varepsilon]\right) \cap [a_i+\varepsilon,b_i-\varepsilon] = \emptyset.
\end{equation}
}
\begin{proof}
We proceed by induction on the number $r$ of intervals.
The case $r=1$ follows from Lemma~\ref{thm:build-large}. 
We now assume $r>1$ and define the following family:
\[
\mathcal{L}=\left\{h \in H \; \Big \vert \;h\left(\bigcup_{i=1}^{r-1}[a_i+\varepsilon,b_i-\varepsilon]\right)\cap \bigcup_{i=1}^{r-1}[a_i+\varepsilon,b_i-\varepsilon] = \emptyset \right\}.
\]
By the induction hypothesis the family $\mathcal{L}$ is non-empty.  We
also note in passing that the set $\mathcal{L}$ is closed under the
operation of passing to inverses.  We will now prove that there is an
element $w_\varepsilon$ in $\mathcal{L}$ with
$w_\varepsilon\left([a_r+\varepsilon,b_r-\varepsilon]\right)\cap
[a_r+\varepsilon,b_r-\varepsilon]=\emptyset$.

For ease of discussion, we denote the orbital $(a_r,b_r)$ by $A_r$.
Let $f \in \mathcal{L}$, if 
$[a_r+\varepsilon,b_r-\varepsilon]\subset \Supp(f)$ then there is some power $n$ so that by setting
$w_\varepsilon = f^n$ we will have found the element we desire, thus,
we assume below that 
$[a_r+\varepsilon,b_r-\varepsilon]\not\subset \Supp(f)$.

Define $\Gamma = \overline{\Supp(f) \cap A_r}$.
There are three possible cases:
\begin{enumerate}
\item Neither $a_r$ nor $b_r$ are in $\Gamma$,
\item Exactly one of $a_r$ and $b_r$ is in $\Gamma$,
\item Both $a_r$ and $b_r$ are in $\Gamma$.
\end{enumerate}
Throughout the cases below we will repeatedly construct a $g\in H$
which will always have an orbital $(s,t)$ containing
$[a_r+\varepsilon,b_r-\varepsilon]$ by evoking
Lemma~\ref{thm:build-large}. We will specify other properties for $g$ as
required by the various cases.

\begin{itemize}
\item[\underline{Case $1$}:]
Possibly by inverting $g$ we can assume that $g$ is increasing on
$(s,t)$, and also by Lemma~\ref{thm:build-large} we can assume that
$s$ is to the left of $\Gamma$ and $t$ is to the right of $\Gamma$
(hence both $s$ and $t$ are fixed by $f$).  Note that for any integers
$m$ and $K$ and for all sufficiently large $n$, the product $f^mg^n
f^{-m} \cdot f^{-K}$ has orbital $(s,t)$ and sends
$[a_r+\varepsilon,b_r-\varepsilon]$ to the right of $b_r-\epsilon$.

\item[\underline{Case $2$}:]
We initially assume $a_r \in \Gamma$. There are two possible subcases.
\begin{enumerate}
\item[(a)] There is an orbital $(a_r,x)$ of $f$, or
\item[(b)] $a_r$ is the accumulation point of a decreasing
  sequence of left endpoints ${x_j}$ of orbitals of $f$ in
  $(a_r,b_r)$.
\end{enumerate}

In case (2.a), possibly by replacing $f$ by its inverse, we can assume
that $f$ is decreasing on the orbital $(a_r,x)$ with $x < b_r$. By
Lemma~\ref{thm:build-large} we can choose $g$ so that $s \in [a_r,x)$
  with $s<a_r+\varepsilon$, $t$ is to the right of $\Gamma$, and $g$
  is increasing on its orbital $(s,t)$ (by inverting $g$ if
  necessary).  For any positive integer $M$ and for all $m \ge M$ we
  have that $f^mg^n f^{-m}$ is increasing on its orbital $(f^m(s),t)
  \supsetneq (s,t)\supsetneq [a_r+\varepsilon,b_r-\varepsilon]$.  It
  is now immediate that for any positive integers $m \ge M$ and $K$
  and for all sufficiently large $n$, the product $f^mg^n f^{-m} \cdot
  f^{-K}$ moves $[a_r+\varepsilon,b_r-\varepsilon]$ entirely off of
  itself to the right.

In case (2.b) we choose an element $x_j$ of the sequence $\{x_p\}$
such that $a_r <x_j <a_r+\varepsilon$.  Moreover, we can choose $g$
increasing so that $a_r<s<x_j$ and $t$ is to the right of $\Gamma$.
For any positive integer $K$ the power $f^{-K}$ fixes the interval
$[x_j,\sup \Gamma] \supseteq [a_r+\varepsilon,b_r-\varepsilon]$
setwise. For any $M>0$ and for any $m\ge M$ the conjugate $f^mg^n
f^{-m}$ has orbital $(f^m(s),t) \supset [x_j,\sup \Gamma]$.
Therefore, there exists an $N>0$ so that for all $n\ge N$ the product
$f^mg^n f^{-m} \cdot f^{-K}$ throws the interval
$[a_r+\varepsilon,b_r-\varepsilon]$ off itself to the right.

If instead in Case 2 we have that $b_r$ is the only endpoint contained
in $\Gamma$ similar arguments prove the existence of a suitable
product $f^mg^n f^{-m} \cdot f^{-K}$ which moves the interval
$[a_r+\varepsilon,b_r-\varepsilon]$ leftward entirely off of itself.

\item[\underline{Case $3$}:]
We have two subcases.

\begin{enumerate}
\item[(a)] $f$ has orbitals $(a_r,x)$ and $(y,b_r)$ with $x<y$, or
\item[(b)] at least one of $a_r$ or $b_r$ is the accumulation point of
a monotone sequence sequence of endpoints ${x_j}$ of orbitals of $f$ in $(a_r,b_r)$, or
\end{enumerate}

In case (3.a) we have that $f$ has orbitals $(a_r,x)$ and $(y,b_r)$
with $x<y$ (if $f$ has $(a_r,b_r)$ as an orbital, then there is a
positive integer $m$ such that $w_\varepsilon:=f^m$ will satisfy our
statement).  We construct $g$ so that it has an orbital $(s,t)$ upon
which it is increasing and where $s \in [a_r,x)$ and $t \in (y,b_r]$.
Possibly by replacing $f$ with its inverse, we can assume that $f$ is
decreasing on the orbital $(a_r,x)$.  We now have two subcases
depending on whether $f$ is increasing or decreasing on $(y,b_r)$.

If $f$ is increasing on $(y,b_r)$, then for any positive integer $M$
and for all $m \ge M$ the conjugate $f^m g f^{-m}$ will have an
orbital containing $(s,t)$.  Given any $K>0$ we can choose an positive
integer $N$ large enough so that, for all $n \ge N$, the element
$f^mg^n f^{-m}$ moves both $x$ and $a_r+\varepsilon$ to the right of
$b_r-\varepsilon$. Under these conditions, the product $f^mg^n f^{-m}
\cdot f^{-K}$ will move $a_r+\varepsilon$ leftward past
$b_r-\varepsilon$.

Assume now that $f$ is decreasing on $(y,b_r)$. There exists an
integer $j>0$ such that $g^j(x)>y$ and so the support of the function
$f^{(g^j)}$ contains the interval $(a_r,y]$.  If $J$ is the orbital of
$f^{(g^j)}$ containing $a_r$, then $J \cup (y,b_r) = (a_r,b_r)$ and so
there exist two positive integers $k_1$ and $k_2$ such that the
support of the function $g^\ast:=(f^{(g^j)})^{k_1} f^{k_2}$ contains
the interval $(a_r,b_r)$.  For any positive integer $M$ and for all $m
\ge M$ the support of $f^m(g^\ast) f^{-m}$ contains $(a_r,b_r)$, hence
for any $K>0$ we can select an integer $n \ge N$ large enough so that
the product $f^m(g^\ast)^n f^{-m} \cdot f^{-K}$ moves the interval
$[a_r+\varepsilon,b_r-\varepsilon]$ off itself.

In case (3.b) we initially assume that $a_r$ is the accumulation point
of a decreasing sequence of left endpoints ${x_j}$ of orbitals of $f$
in $(a_r,b_r)$. Now, either $f$ has a fixed point $y \ge
b_r-\varepsilon$ or it has an orbital $(y,b_r)$ with
$y<b_r-\varepsilon$.  In the second case we will assume $f$ is
increasing on its orbital $(y,b_r)$ (possibly by replacing $f$ by its
inverse).  In either case we choose $g$ decreasing on $(s,t)$ so that
$t > y$ and $t>b_r-\varepsilon$.  We also assume $g$ is chosen so that
$s$ is to the left of a fixed point of $f$ which is to the left of
$a_r+\varepsilon$.  Now by our choices it is easy to see that given
any positive $M$ and $m>M$ and any positive $K$ we have
\begin{enumerate}
\item $f^{-K}(b_r-\varepsilon)<f^m(t)$,
\item $f^{-K}(a_r+\varepsilon)>f^m(s)$, and
\item there is positive $N$ so that for all $n> N$ we have $f^mg^nf^{-m}\cdot f^{-K}(b_r-\varepsilon)<a_r+\varepsilon$.
\end{enumerate}

A similar (reflected) argument can be made if $b_r$ is the
accumulation point of an increasing sequence of right endpoints
${x_j}$ of orbitals of $f$ in $(a_r,b_r)$.

\end{itemize}

By Lemma~\ref{thm:product-power-conjugate} there exists an $M_0$ such
that for all $m \ge M_0$ we can find a $K_0>0$ such that for all $k
\ge K_0$ we can find an $N_0>0$ so that for all $n \ge N_0$ the
product $f^m g^n f^{-m} \cdot f^{-k}$ has support containing
$\bigcup_{i=1}^{r-1}[a_i+\varepsilon,b_i-\varepsilon]$. By the
analysis in this proof we know we can choose an $M \ge M_0$ such that
for any $m \ge M$ we can find a $K \ge K_0$ and $N \ge N_0$ (depending
on $K$) so that for all $n \ge N$ the product $w_\varepsilon:=f^m g^n
f^{-m} \cdot f^{-K}$ throws $[a_r+\varepsilon,b_r-\varepsilon]$
entirely off of itself.
\end{proof}

We are finally in position to prove the Lemma~\ref{thm:FIP} from
our introduction.
\begin{proof}[Proof of Lemma~\ref{thm:FIP}]
We argue via induction on $n$, with the case $n=2$ being true by Lemma~\ref{thm:free-group}.
Let $g_1, \ldots g_n \in G_0$ and define $H:=\langle g_1, \ldots,
g_{n-1} \rangle$.

Write $S^1 \setminus \Fix(H)=\bigcup I_\alpha$ and $S^1
\setminus \Fix(g_n)=\bigcup J_\beta$, for suitable families
of open intervals $\{I_\alpha\},\{J_\beta\}$.

We assume, by contradiction, that $\Fix(H) \cap
\Fix(g_n) = \emptyset$, hence we have $S^1 \subseteq
\left(\bigcup I_\alpha\right) \cup \left(\bigcup J_\beta\right).$
By the compactness of $S^1$ and there are indices $r$ and $s$ so that we
can write $S^1=I_1 \cup \ldots \cup I_r \cup J_1 \cup \ldots \cup
J_s$.

Let $I_i=(a_i,b_i)$ and notice that $\Fix(H) \cap \left(
\bigcup_{i=1}^r I_i \right) = \emptyset$, so we can apply Lemma
\ref{throw-off} to build an element $w_\varepsilon \in H$ such that
$\bigcup_{i=1}^r (a_i+\varepsilon,b_i-\varepsilon) \subseteq
\Supp(w_\varepsilon).$ We choose $\varepsilon>0$ to be
small enough so that $\Fix(g_n) \subseteq \bigcup_{i=1}^r
(a_i+\varepsilon,b_i-\varepsilon)$ thus implying
$\Fix(w_\varepsilon)\cap \Fix(g_n)=\emptyset$. By
Lemma \ref{thm:free-group} we can find a non-abelian free group inside
$\langle w_\varepsilon,g_n\rangle$, contradicting the assumption on
$G$.
\end{proof}

By compactness of $S^1$, the previous lemma immediately implies:

\corollary{\label{thm:compactness-global} Let $G \le \homeooS$ with no
  non-abelian free subgroups.  Then 
\begin{enumerate}
\item $G_0$ admits a global fixed point,
  i.e., $\Fix(G_0)\ne \emptyset$, and so
\item $G_0$ is a normal subgroup
  of $G$.
\end{enumerate}  }

Another application of the compactness is:
\claim{Let $f \in \homeooS$, then for any $0<\varepsilon<1$
there exists integer $n >0$ and a point $x \in S^1$ such that
the distance between $x$ and $f^n(x)$ is less than $\varepsilon$, i.e.,
$\wh{f}^n(\wh{x}) = \wh{x} + k + \delta$ for some integer $k$ and
$|\delta| <\varepsilon$.
}
\begin{proof}
Let $y$ be any point on $S^1$. The sequence
$\{f^n(y)\}_n$ contains a converging subsequence
$\{f^{n_i}(y)\}_i$. Therefore there exist $i <j$ such that
distance between $f^{n_i}(y)$ and $f^{n_j}(y)$ is less the $\varepsilon$.
Thus, we can take $x:=f^{n_i}(y)$ and $n= n_j-n_i$.
\end{proof}

\lemma{
\label{rationalInsertion}
Given $f,g \in \homeooS$ such that $\wh{f} < \wh{g}$, then there
exists a function $h\in \homeooS$ with rational rotation number and
such that $\wh{f}<\wh{h}<\wh{g}$.}

\begin{proof}
%
Let $\varepsilon$ be the minimal distance between $\wh{f}$ and $\wh{g}$, i.e.,
\[
\varepsilon = \frac{1}{2}\min_{t \in [0,1]}\left\{|\wh{f}(t)-\wh{g}(t)|\right\}
\]
and let
$\wh{h_0}:=(\wh{f}+\wh{g})/2$.
Choose $x$ and $n$ be the ones given by the claim for the function $h_0$ and the value $\varepsilon/3>0$, i.e.,
$|\wh{h_0}^n(\wh{x}) - \wh{x} - k | < \varepsilon/3$ for some integer $k$.
Consider the family of functions $\wh{h_t}(s) := \wh{h_0}(s) + t$
and their powers $\wh{h_t}^n$. The monotonicity of $\wh{h_t}$
implies that for any $t>0$, we have
$$
\wh{h_t}^n(s) = \wh{h_t}(\wh{h_t}^{n-1}(s)) =
\wh{h_0}(\wh{h_t}^{n-1}(s)) +t \geq
\wh{h_0}(\wh{h_0}^{n-1}(s)) =
\wh{h_0}^n(s)+t.
$$
Similarly we have $\wh{h_t}^n(s) \leq \wh{h_0}^n(s) +t$ if $t <0$.
The intermediate value theorem applied to the function $t \to \wh{h_t}^n(\wh{x})$
implies that there exists a $t$ such that $|t| \leq \varepsilon/3$ and
$\wh{h_t}^n(\wh{x}) - \wh{x}=k$ is an integer, i.e., $x$ is a periodic point for $h_t$.
Hence $h_t$ has rational rotation number.
By construction $\wh{h_t}$ is very close to $\wh{h_0}$, therefore it is between $\wh{f}$ and $\wh{g}$.
\end{proof}

The proof of Lemma~\ref{thm:RotHom} involves observing that the
element $(fg)^n$ can be rewritten $f^n g^n h_n$ for some suitable
product of commutators $h_n \in [G,G]$; if we prove that $[G,G]$ has a
global fixed point $s$ we can compute the rotation number on $s$, so
that $(fg)^n(s)=(f^n g^n)(s)$. The next lemma, together with Corollary
\ref{thm:compactness-global}, shows that this is indeed the case.
\lemma{
\label{thm:same-rot-fix-point}
Let $G \le \homeooS$ and let $f,g \in G$. Suppose one
  of the following two cases is true:
\begin{enumerate}
\item $G$ has no non-abelian free subgroups and $\rot(f)=\rot(g) \in \mathbb{Q}/\mathbb{Z}$, or
\item $\rot(f)=\rot(g) \not \in \mathbb{Q} / \mathbb{Z}$.
\end{enumerate}
Then $fg^{-1} \in G_0$. 
}

\begin{proof}
(1) Assume $\rot(f)=\rot(g)=p/q \in \mathbb{Q} / \mathbb{Z}$ with
  $p,q$ positive integers
 and that $G$ has no non-abelian free
  subgroups.

In this case, $f^q$ and $g^q$ have fixed points in $S^1$.  Now,
$\wh{f}^q(\wh{x})=\wh{x}+p$ and $\wh{g}^q(\wh{y})=\wh{y}+p$ for any $x
\in \Fix(f^q)$ and $y \in \Fix(g^q)$, by Lemma~\ref{thm:lift-control}(4).
 In particular, $f^q$ and $g^q$ must have
a common fixed point $s \in S^1$ by Lemma~\ref{thm:free-group} (in the
case that one of $f^q$ or $g^q$ is the identity map, then it is
immediate that $f^q$ and $g^q$ have a common fixed point) and then for
this $s$ we must have $\wh{f}^q(\wh{s}) = \wh{s} +p = \wh{g}^q(\wh{s})$.

  Suppose now that $fg^{-1} \not \in G_0$. In this case, either
$\wh{f}>\wh{g}$ or $\wh{f}<\wh{g}$.  We suppose without meaningful
loss of generality that the latter is true.  However, $f<q$ implies $\wh{f}^q<\wh{g}^q$,
which is impossible as $\wh{f}^q(\wh{s})=\wh{s}+p=\wh{g}^q(\wh{s})$.

(2) Assume now that $\rot(f)=\rot(g) \not \in \mathbb{Q} /
  \mathbb{Z}$.

Suppose $fg^{-1}\not \in G_0$.  Again, either $\wh{f}<\wh{g}$ or $\wh{g}
<\wh{f}$.  Without meaningful loss of generality we suppose that
$\wh{f}<\wh{g}$.  By Lemma~\ref{rationalInsertion} we can find a map
$h\in \homeooS$ with $\wh{f}<\wh{h}<\wh{g}$ where $h$ has
rational rotation number.  However, this is impossible since
$\wh{f}<\wh{h}<\wh{g}$ guarantees us that $\rot(f) \leq \rot(h) \leq
\rot(g)=\rot(f)$, so that all three rotation numbers must be equal.

In both (1) and (2), we ruled out the possibility that $fg^{-1}\not\in
G_0$, thus we must have that $fg^{-1}\in G_0$.
\end{proof}
\corollary{
\label{thm:commutator-fixed}
Let $G \le \homeooS$ with no non-abelian free subgroups, then we
have $[G,G] \le G_0$.
}

The following Lemma is an easy consequence of the definition of
lift of a map and Corollary~\ref{thm:free-group-2} and we omit its
proof (it can be found in~\cite{MatucciThesis}).

\lemma{
\label{thm:lift-commutator-fixed}
Let $G \le \homeooS$ with no non-abelian free subgroups.
  Let $u,v \in G$ and $s \in S^1$ be a fixed point of
  \mk{the commutator} $[u,v]$.  Then
  $\wh{s}$ is a fixed point for $[U,V]$, for any $U$ lift of $u$ and
  $V$ lift of $v$ in $\Homeoo(\mathbb{R})$.
}

We are now ready to give a proof the main result of this section.

\begin{proof}[Proof of Lemma~\ref{thm:RotHom}] Let $f,g \in G$. We write the power
$(fg)^n=f^n g^n h_n$ where $h_n$ is a suitable product of commutators \mk{(involving $f$ and $g$)}
  used to shift the $f$'s and $g$'s leftward.  Since $h_n \in [G,G]
  \le G_0$ for all positive integers $n$ then, if $s \in S^1$ is a
  global fixed point for $G_0$, we have $h_n(s)=s$. Similarly, we
  observe that $(\wh{f}\wh{g})^n=\wh{f}^n \, \wh{g}^n \, H_n$ where
  $H_n$ is a suitable product of commutators and $H_n$ is a lift for
  $h_n$. By Lemma~\ref{thm:lift-commutator-fixed} we must have that
  $H_n(\wh{s})=\wh{s}$ for all positive integers $n$. Thus we observe
  that:
\[
(\wh{f}\wh{g})^n(\wh{s}) = \wh{f}^n \, \wh{g}^n \, H_n(\wh{s}) =
\wh{f}^n \, \wh{g}^n (\wh{s}).
\]
We now find upper and lower bounds for $\wh{f}^n \, \wh{g}^n (\wh{s})$. Observe that,
for any two real numbers $a,b$ we have that
\[
\wh{f}^n(a) + b - 1 < \wh{f}^n(a) + \lfloor b \rfloor \le  \wh{f}^n(a+b) <
\wh{f}^n(a)+ \lfloor b \rfloor +1 \le \wh{f}^n(a) + b +1
\]
where $\lfloor \cdot \rfloor$ denotes the floor function. By applying this inequality to
$\wh{f}^n \, \wh{g}^n (\wh{s}) = \wh{f}^n(\wh{s}+(\wh{g}^n(\wh{s})-\wh{s}))$ we get
\begin{eqnarray*}
\wh{f}^n(\wh{s})+\wh{g}^n(\wh{s})-\wh{s} - 1 \le
\wh{f}^n(\wh{s}+(\wh{g}^n(\wh{s})-\wh{s})) \le
\wh{f}^n(\wh{s})+\wh{g}^n(\wh{s})-\wh{s}+1.
\end{eqnarray*}
We divide the previous inequalities by $n$, and get
\[
\frac{\wh{f}^n(\wh{s})+\wh{g}^n(\wh{s})-\wh{s}-1}{n} \le \frac{(\wh{f} \, \wh{g})^n(\wh{s})}{n} \le \frac{\wh{f}^n(\wh{s})+\wh{g}^n(\wh{s})-\wh{s}+1}{n}.
\]
By taking the limit as $n \to \infty$ of the previous expression, we
immediately obtain $\rot(fg)=\rot(f)+\rot(g)$.
\end{proof}

\corollary{
\label{thm:G_0-normal}
Let $G \le \homeooS$ with no non-abelian free subgroups. Then $\rot:G \to \mathbb{R}/\mathbb{Z}$
is a group homomorphism and
\begin{enumerate}
\item $\ker(\rot)=G_0$,
\item  $G/G_0 \cong \rot(G)$.
\item for all $f,g \in G$, $fg^{-1} \in G_0$ if and only if $\rot(f)=\rot(g)$.
\end{enumerate}
}

\section{Applications: Margulis' Theorem \label{sec:margulis}}

In this section we show how the techniques developed in Section~\ref{sec:rot-homo} 
yield two results for groups of homeomorphisms of
the unit circle.  One of these results is Margulis' Theorem (Theorem 1.10) which
states that every group $G$ of orientation-preserving homomorphisms of the unit circle $S^1$
either contains a non-abelian free subgroup or admits a $G$-invariant probability
measure on $S^1$.

\begin{proof}[Proof of Theorem~\ref{MargulisThm}]
We assume that $G$ does not contain free subgroups, so that the $\rot$ map is a group
homomorphism, by Lemma~\ref{thm:RotHom}. The proof divides into two
cases.

\medskip
\noindent \emph{Case 1: $G /G_0$ is finite.}\\
Let $s \in \Fix(G_0)$ and consider the finite orbit $s^G$.
Then for every subset $X \subseteq S^1$ we assign:
\[
\mu(X)=\frac{\# \, s^G \cap X}{\# \, s^G}.
\]
This obviously defines a $G$-invariant probability measure on $S^1$.

\medskip
\noindent \emph{Case 2: $G/G_0$ is infinite and therefore $\rot(G)$ is dense in $\mathbb{R} / \mathbb{Z}$.}\\

Fix $s \in \Fix(G_0)$ as an origin and
\mk{identify $S^1$ with}
$[0,1]$.
 We regard $s^G$ as a subset of $[0,1]$ and define the map
$\varphi:s^G \to \rot(G)$, given by $\varphi(s^g)=\rot(g)$, for any $g
\in G$. It is immediate that
 $\varphi$ is well-defined and
order-preserving on $s^G \subseteq [0,1]$. We take the
\mk{continuous extension} 
of
this map, by defining the function:
\[
\begin{array}{cccc}
\ov{\varphi}: &  [0,1] & \longrightarrow & [0,1] \\
              &    a   & \longmapsto     & \sup \{\rot(g) \; | \; s^g \le a,  \; \; g \in G \}.
\end{array}
\]
By construction, the function $\ov{\varphi}$ is
non-decreasing. Moreover, since the image of $\ov{\varphi}$ contains
$\rot(G)$, it is dense in $[0,1]$. Since $\ov{\varphi}$ is a
non-decreasing function whose image is dense in $[0,1]$,
$\ov{\varphi}$ is a continuous map. This allows us to define the
Lebesgue-Stieltjes
 measure associated to $\ov{\varphi}$ on the Borel
algebra of $S^1$ (see~\cite{malliavin}),
 that is, for every
half-open interval $(a,b] \subseteq S^1$ we define:
\[
\mu((a,b]):=\ov{\varphi}(b)-\ov{\varphi}(a).
\]
Since the $\rot$ map is a homomorphism, it is straightforward to see
that the measure $\mu$ is $G$-invariant. By definition, $\mu(S^1)=1$
and $\mu(\{p\})=0$, for every point $p \in S^1$.
\end{proof}

\mk{\medskip}

Next, we impose a categorical restriction on our group of
homeomorphisms, so that Denjoy's Theorem applies.  Under these
conditions, the  existence of an element with irrational rotation number
yields an analog of the Tits' alternative --- either the group is abelian
or it contains a non-abelian free group.

\medskip

\begin{proof}[Proof of Theorem~\ref{thm:strongTits}]
Let us suppose $G$ contains no non-abelian free subgroups, and let
$s\in \textrm{ Fix(} G_0$).  By Denjoy's Theorem there is a $z$
in $\homeooS$ so that $g^z$ is a pure rotation (by an irrational
number).  Thus, the orbits of $g$ are dense in $S^1$ so in particular
the orbit of $s$ under the action of $g$ is dense in $S^1$.  Since
Fix($G_0$) must be preserved as a set by the action of $G$, we see that
$G_0$ must be the trivial group. By Corollary~\ref{thm:G_0-normal}, we have $G \cong \rot(G) \le \mathbb{R} / \mathbb{Z}$ and that $G$ is contained in
$C_{\homeooS}(g) \cong \mathbb{R} / \mathbb{Z}$.
\end{proof}

\section{Structure and Embedding Theorems}
\label{sec:structure-embedding}
\subsection{Structure Theorems}

We start the section with our main result which classifies the
structure of subgroups of $\homeooS$ with no non-abelian free
subgroups. We consider an orbit $s^G$ of a point $s$ of $\Fix(G_0)$
under the action of $G$ (recall that $\ov{s^G} \subseteq \Fix(G_0)$),
then we choose a fundamental domain $D$ for the action of $G$ on $S^1
\setminus \ov{s^G}$.  Since the subset $S^1 \setminus \ov{s^G}$ is
open, the fundamental domain will be given by a union of intervals. By
restricting $G_0$ to this fundamental domain we get a group $H_0$
which acts as a set of homeomorphisms of a disjoint union of
intervals.  We will prove that if $G \le \homeooS$ without non-abelian
free subgroups then either $G$ is abelian or $G$ can be embedded into
the wreath product $H_0 \wr (G/G_0)$.

\remark{
\label{thm:PL-is-always-rational}
Note that by Theorem~\ref{thm:irrational-abelian} (a consequence of
Denjoy's Theorem), if $G \le \PL_+(S^1)$ is non-abelian with no
non-abelian free subgroups, then $Q$ is isomorphic to a subgroup of
$\mathbb{Q} / \mathbb{Z}$.  }

\begin{proof}[Proof of Theorem~\ref{structureThm}]
If $G_0=\{\mathrm{id}_{S^1}\}$, then Corollary~\ref{thm:G_0-normal}
implies $G \cong G/G_0 \cong \rot(G) \le \mathbb{R}/\mathbb{Z}$.
Now suppose $G_0$ non-trivial, so that $\Fix(G_0) \ne S^1$ and
define $P=G/G_0$. Let $s \in \Fix(G_0)$.  Note that $P$ acts
on Fix($G_0$) and consider the open subset $S^1 \setminus \ov{s^P}$,
where $s^P$ is the orbit of $s$ under the action of $P$.  The set $S^1
\setminus \ov{s^P}$ is a collection of at most countably many disjoint
open intervals.  We observe that $P$ also acts on $S^1 \setminus
\ov{s^P}$ thought of as a set whose elements are open intervals.  We
can define a fundamental domain for the action of $P$ on $S^1
\setminus \ov{s^P}$ as the union $D=\bigcup_{i \in \mathfrak{N}} I_i$
of a collection $\{I_i\}_{i \in \mathfrak{N}}$ of at most countably
many intervals $I_i$ such that
\begin{eqnarray*}
k_1(D) \cap k_2(D) = \emptyset, \, \, k_1 \ne k_2,  \\
\\
S^1 \setminus \ov{s^P} = \bigcup_{k \in P} k(D).
\end{eqnarray*}

\medskip
We give proof of Claim 
\ref{thm:end-proof} below, and leave the remaining claims to the reader.

\medskip 
\claim{
The fundamental domain $D$ exists.
\label{thm:claim1}
}

\noindent Since $\ov{s^P} \subseteq \Fix(G_0)$ we have
\[
S^1 \setminus \bigcup_{k \in P} k(D) \subseteq \Fix(G_0).
\]

\medskip
\claim{Define $H_0 \le \homeooS$ to be the subgroup generated by
  functions $f$ such that there exists a function $g_f \in G_0$ with
  $f$ the restriction of $g_f$ on $D$ and the identity on $S^1
  \setminus D$.  Then $H_0 \hookrightarrow \prod_{i\in \mathfrak{N}}
  \Homeoo(I_i)$, since $D = \bigcup_{i\in \mathfrak{N}} I_i$.
  Similarly for every $k \in G/G_0$, there is an embedding $H_0
  \hookrightarrow \prod_{i\in \mathfrak{N}} \Homeoo(k^{-1}(I_i))$.  }
\remark{We will call the image group of this last embedding $H_0^k$.}


It is important to notice that $H_0$ is not necessarily contained in $G_0$,
since $H_0$ has its
support in $D$, while an element of $G_0$ has support in $\bigcup_{k \in P} k(D)$.

\medskip
\claim{The conjugates of $H_0$ under $P$ commute, and the group
  $\widetilde{H}:=\langle H_0^s \mid s\in G \rangle \simeq
  \bigoplus_{k \in P} H_0^k$ is normalized by $G$.  Moreover, the
  group $H:=\prod_{k \in P} H_0^k$, thought of as a subgroup in
  $\Homeoo(S^1)$, contains $\widetilde{H}$ and is also normalized by $G$.
\label{thm:build-H}}

%

\medskip
We define the following subgroup
\[
E:=\langle G, H \rangle \le \homeooS
\]
and observe that, since $G$ normalizes $H$ by Claim \ref{thm:build-H},
the group $H$ is normal in $E$ and we have the following exact
sequence:
\[
1 \to H \overset{i}{\to} E \overset{\pi}{\to} E/H \to 1
\]
where $i$ is the inclusion map and $\pi$ is the natural projection
$\pi:E \to E/H$.  Notice that $E/H \cong G/(G\cap H)$ and $G \cap H =
G_0$, by definition of $G_0$.  Thus, $E/H \cong G/G_0=P$ , so we can rewrite the sequence as
\[
\qquad\qquad\qquad\qquad\qquad\qquad 1 \to H \overset{i}{\to} E \overset{\pi}{\to} P \to 1.\qquad\qquad\qquad\qquad\qquad\left(*\right)
\]

Since $G$ is a subgroup of $E$, the conclusion of the theorem will
follow if we can show that $E\cong H_0\wr P$, where $H$ in the exact
sequence $\left(*\right)$ above plays the role of the base group.  In
this case, the semi-direct product structure of $E$ enables us to find
a splitting $\phi:P\to E$ of the exact sequence $\left(*\right)$ so
that if we set $Q = \mathrm{Im}(\phi) \cong P$ we will have the remaining
points of our statement.

\medskip
\claim{The group $H\rtimes P\cong H_0\wr P$ is the only
  extension of $\prod H_0^k = H$ by $P$, where $P$ acts on $H$ by permuting the copies of $H_0$.
\label{thm:end-proof}}
\begin{proof}  
By a standard result in cohomology of groups (see Theorem~11.4.10 in
\cite{rob}), if we can prove that $H^2(P,Z(\prod H_0^k))=0$ (where
$Z(\prod H_0^k)$ denotes the center of $\prod H_0^k$), there can be
only one possible extension of $\prod H_0^k$ by $P$.  We observe that
$H \rtimes P \simeq H_0 \wr P$ is one such extension, so it suffices to prove that
$H^2(P,Z(\prod H_0^k))=0$. We use Shapiro's Lemma to compute this
cohomology group (see Proposition 6.2 in \cite{brown1}). We have
\begin{eqnarray*}
H^2(P,Z(\prod H_0^k))=H^2(P,\prod Z(H_0)^k)= \\
=H^2(P,\mbox{Coind}^P_{\{e\}} Z(H_0)) = H^2(\{e\},Z(H_0))=0,
\end{eqnarray*}
which completes the proof of the claim.
\end{proof}

\medskip 
\noindent \emph{End of the proof of Theorem~\ref{structureThm}.}
Since $E\cong H_0\wr P$, there is a splitting $\phi:P\to
E$ of the exact sequence $\left(*\right)$ so that $E=\langle
H,Q\rangle\cong H\rtimes Q$ where $Q =$ Im($\phi$) $\cong P$.

\end{proof}

\remark{We observe that the wreath product in the previous result is
  unrestricted; the elements of $\homeooS$ can have infinitely many
  ``bumps'' and so the elements of $G_0$ can be non-trivial on
  infinitely many intervals. On the other hand, if we assume $G \le
  \PL_+(S^1)$, this would imply that any element in $G_0$ is
  non-trivial only at finitely many intervals, and so $G_0$ can be
  embedded in the direct sum $\bigoplus$.  This argument explains why
  the wreath product in Theorem~\ref{thm:first-embedding} is
  unrestricted whereas the ones in
  Theorems~\ref{thm:embedding-thompson} and~\ref{thm:embedding-PL} are
  restricted.}

We now obtain structure results about solvable subgroups of
$\PL_+(S^1)$.
 Following the first author in~\cite{bpasc}, we define
inductively the following family of groups. Let $G_0=1$
 and, for $n
\in \mathbb{Z}_+$, we define $G_n$ as the direct sum of infinitely
many copies of the group $G_{n-1} \wr \mathbb{Z}$:
\[
G_n := \bigoplus_{d \in \mathbb{Z}} \left(G_{n-1} \wr \mathbb{Z} \right).
\]
We recall the following classification.

\theoremname{Bleak~\cite{bpasc}} 
{ 

Let $H$ be a solvable group with derived length $n$. Then, $H$ embeds
in $\ploi$ if and only if $H$ embeds in $G_n$.

}

Using Theorem~\ref{structureThm}
and Remark~\ref{thm:PL-is-always-rational}, we are able to extend this result \mk{to }
obtain Theorem 1.3 from the introduction.

\medskip

There is also a non-solvability criterion for subgroups of $\PL_+([0,1])$.
Let $W_0=1$
 and, for $n \in \mathbb{N}$, we define $W_i= W_{i-1}
\wr \mathbb{Z}$.
 We build the group
\[
W := \bigoplus_{i \in \mathbb{Z}} W_i.
\]

The following is the non-solvability criterion mentioned above.

\theoremname{Bleak~\cite{bpnsc}}
{Let $H \le \PL_+([0,1])$. Then $H$ is non-solvable if and only if it
contains a subgroup isomorphic to $W$.}

Using this result and Theorem 1.1, one immediately derives a Tits'
alternative for subgroups of $\PL_+(S^1)$; Theorem 1.4 from the
introduction.

\medskip

\subsection{Embedding Theorems}

We now turn to prove existence results and show that subgroups with wreath product structure do exist
in $\homeooS$ and in $\PL_+(S^1)$.

\medskip

\remark{The same result is true for any
$H_0$ that can be embedded in $\prod \Homeoo(I_i)$ (following the notation of Theorem~\ref{structureThm})
and our proof can be extended without much effort,
however we prefer to simplify the hypothesis in order to keep the proof cleaner.
\mk{Alternatively, we can use the existence of embedding $\prod_{i\in K} \Homeoo(I_i) \to \homeooI$
if $K$ is countable.}
}

\begin{proof}[Proof of Theorem 1.5] We divide the proof into two cases: $K$ infinite and $K$ finite.
If $K$ is infinite, we enumerate the elements of $K=\{k_1, \ldots, k_n, \ldots \}$ and we
choose the sequence:
\[
\frac{1}{2}, \frac{1}{2^2}, \ldots, \frac{1}{2^n}, \ldots
\]
We identify $S^1$ with the interval $[0,1]$ to fix an origin and an orientation of the unit circle.
$K$ is countable subgroup of $\mathbb{R}/\mathbb{Z}$, so it is non-discrete and therefore it is dense in $S^1$.
Now define the following map:
\[
\begin{array}{cccc}
\varphi: & [0,1] =S^1 & \longrightarrow & [0,1] =S^1\\
         & x          & \longmapsto     & \sum_{k_i<x} \frac{1}{2^i}
\end{array}
\]
(where $k_i < x$ is written with respect to the order in $[0,1]$). It is immediate from
the definition to see that the map is order-preserving and it is injective, when restricted
to $K$.

For small enough $\varepsilon>0$ we have
\[
\varphi(k_1 + \varepsilon) = \sum_{k_i< k_1 + \varepsilon} \frac{1}{2^i}.
\]
If we let $\varepsilon \to 0$, we then see that
\begin{eqnarray*}
\alpha : = \varphi(k_1) < \varphi(k_1 + \varepsilon) \underset{\varepsilon \to 0}{\longrightarrow}
\sum_{k_i \le k_1} \frac{1}{2^i} = \alpha + \frac{1}{2}.
\end{eqnarray*}
But now, as $\varphi$ is non-decreasing, we must have $(\alpha, \alpha + \frac{1}{2}) \cap \varphi(K)=\emptyset$.
More generally, it follows that:
\[
\bigcup_{i \in \mathbb{N}} \left(\varphi(k_i),\varphi(k_i)+\frac{1}{2^i}\right) \cap \ov{\varphi(K)} = \emptyset
\]
\claim{The unit circle can be written as the disjoint union
\[
S^1=\bigcup_{i \in \mathbb{N}} \left(\varphi(k_i),\varphi(k_i)+\frac{1}{2^i}\right) \cup \ov{\varphi(K)}.
\]
}
\begin{proof}
 Let $A := \bigcup_{i \in \mathbb{N}} \left(\varphi(k_i),\varphi(k_i)+\frac{1}{2^i}\right)$
and let $x_0 \not \in A$. Let $\varepsilon >0$ be given.  We want to
prove that we have $\varphi(K)\cap (x_0-\varepsilon,x_0+\varepsilon)
\ne \emptyset$.

Suppose
$(x_0-\varepsilon,x_0+\varepsilon)\cap A=\emptyset$, then we have
\[
1=m([0,1]) \ge m((x_0-\varepsilon,x_0+\varepsilon)) + m(A)=
2\varepsilon + \sum_{i=1}^\infty \frac{1}{2^i} =2\varepsilon +1 > 1
\]
where $m$ is the Lebesgue measure on $[0,1]$.  In particular, we must
have that \mbox{$(x_0-\varepsilon,x_0+\varepsilon)\cap A$} is not empty.

From the above, we know there is an index $i$ with $k_i \in K$ so that
\[
(x_0-\varepsilon, x_0+\varepsilon) \cap \left(\varphi(k_i),\varphi(k_i) + \frac{1}{2^i} \right) \ne \emptyset.
\]

There are three cases of interest.
\begin{enumerate}
\item [(a)] $\varphi(k_i) \in (x_0-\varepsilon,x_0+\varepsilon)$.

In this case, as $\varepsilon >0$ was arbitrary, we have shown that
$x_0$ is in the closure of $\varphi(K)$.

\item [(b)] $\varphi(k_i) + \frac{1}{2^i} \in
  (x_0-\varepsilon,x_0+\varepsilon)$.

Let $\{k_{i_r}\} \subseteq K \subseteq [0,1]$ be a decreasing sequence
converging to $k_i$. Then, $\lim_{r \to \infty}
\varphi(k_{i_r})=\varphi(k_i)+\frac{1}{2^i}$ and so there is an $r$ such
that $\varphi(k_{i_r}) \in (x_0-\varepsilon,x_0+\varepsilon)$,
returning us to the previous case.

\item [(c)] $(x_0-\varepsilon,x_0+\varepsilon) \subseteq
  \left(\varphi(k_i),\varphi(k_i)+\frac{1}{2^i}\right)$.

This implies that $x_0 \in
\left(\varphi(k_i),\varphi(k_i)+\frac{1}{2^i}\right) \subseteq A$,
which contradicts our definition of $x_0$, so this case cannot occur.

\end{enumerate}

In all possible cases above, we have that $x_0$ is in the closure of
$\varphi(K)$, so our claim is proven.
\end{proof}

We can visualize the set $C:=\ov{\varphi(K)}$ as a Cantor set.  If we
regard $[0,1]$ as $S^1$, then the group $K$ acts on $S^1$ by rotations
and so each $k \in K$ induces a map $k:C \to C$. Now we extend this
map to a map $k:S^1 \to S^1$ by sending an interval
$X_i:=\left[\varphi(k_i),\varphi(k_i)+\frac{1}{2^i}\right] \subseteq
S^1 \setminus C$ linearly onto the interval
$k(X_i):=\left[\varphi(k_j),\varphi(k_j)+\frac{1}{2^j} \right]$, where
$k_j=k+k_i$ according to the enumeration of $K$. Thus we can
identify $K$ as a subgroup of $\homeooS$.

 We squeeze the interval $I$ into $X_1$ and regard the group
 $H_0$ as a subgroup of $\{g \in \homeooS \; | \; g(x)=x, \forall
 x \not \in X_1\} \cong \Homeoo(X_1)$ (we still call $H_0$ this
 subgroup of $\homeooS$).

We now consider the subgroup $H\le\homeooS$ whose elements are
fixed away from all conjugates of $X_1$ (by the action of $K$), and
restrict to elements of $H_0^k$ over $k(X_1)$.  Thus, $H$ is the group
we obtain spreading the action of $H_0$ over the circle through
conjugation by elements of $K$ (where these elements are allowed to be
non-trivial even across infinitely many such conjugate intervals).

Since supp$(H_0^k) \subseteq k(X_1)$ for any $k \in K$, the groups
$H_0^k$ have disjoint support hence they commute pairwise
thus $H\cong
 \prod_{k \in K}H_0^k$. Moreover, the conjugation
action of $K$ on $H$ permutes the
subgroups $H_0^k$. If follows that
\[
\langle H,K\rangle= H_0 \wr K \hookrightarrow \homeooS.
\]
In case $K=\{k_1,\ldots,k_n\}$ is finite, then it is a closed subset
of $S^1$. We define $X_i:=(k_i,k_{i+1})$, for $i=1,\ldots,n$, where
$k_{n+1}:=k_1$.  We can copy the procedure of the infinite case, by
noticing that $S^1=\bigcup_{i=1}^n X_i \cup K$ and embedding $H_0$
into subgroups of $\homeooS$ isomorphic with $\Homeoo(X_i)$.
\end{proof}

We now follow the previous proof, but we need to be more careful in
order to embed Thompson's group $T$ into $\PL_+(S^1)$ (see Section
\ref{sec:back-tools} for the definition
 of Thompson's groups $T$ and
$F$).

 \proposition{There is an embedding
  $\varphi:\mathbb{Q}/\mathbb{Z} \hookrightarrow T$ such
 that
  $\rot(\varphi(x))=x$ for every $x \in \mathbb{Q}/\mathbb{Z}$ and
  there is an interval $I \subseteq S^1$
 with dyadic endpoints such
  that $\varphi(x)I$ and $\varphi(y)I$ are disjoint, for all
 $x,y
  \in \mathbb{Q}/\mathbb{Z}$ with $x \ne y$.}

\begin{proof} \emph{Outline of the idea.}
We consider the set of elements $\{x_n=1/n! \mid n \in \mathbb{N}\}$ of $\mathbb{Q}$ which are the primitive
$n!$-th roots of $1$ in $\mathbb{Q}$
with respect to addition. That is, $n x_n = x_{n-1}$ for each $n$. We want to send each $x_n$ to a homeomorphism $X_n$
of $T$ with $\rot(X_n)=1/n!$ and such that $X_n^n = X_{n-1}$ and $(X_n)^{n!}=\mathrm{id}_{S^1}$.  Then, as
$\langle x_n \mid n \in \mathbb{N}\rangle = \mathbb{Q}/\mathbb{Z}$, we will have an embedding
$\mathbb{Q}/\mathbb{Z} \hookrightarrow T$.

\medskip
\noindent \emph{Notation for the proof.} For every positive integer
$n$ we choose and fix a partition $P_n$ of the unit interval $[0,1]$
into $2n-1$ intervals whose lengths are all powers of $2$. To set up
notation, we always assume we are looking at $S^1$ from the origin of
the axes: from this point of view right will mean clockwise and left
will mean counterclockwise and we will always read intervals
clockwise. We are now going to use the partitions $P_n$ of the unit
interval to get new partitions of the unit circle.

Assume we have a partition of $S^1$ in $2m$ intervals, we define a
``shift by 2'' in $T$ to be the homeomorphism $X$ which permutes the
intervals of the partition cyclically such that $\rot(X)=1/m$
and $X^m=\mathrm{id}_{S^1}$. In other words, ``shift by 2'' sends
an interval $V$ of the partition linearly to another interval $W$
which is 2 intervals to the right of $V$.

\medskip
\noindent \emph{Defining the maps $X_n$.}  We want to build a sequence
of maps $\{X_n\}$ each of which acts on a partition of $S^1$ consisting of
$2(n!)$ intervals $J_{n,1},I_{n,1} \ldots, J_{n,n!},I_{n,n!}$ ordered so that each is to the right of the previous.  The map
$X_n$ will act as the ``shift by 2'' map on this partition.  We define
$X_1=\mathrm{id}_{S^1}$.  To build $X_2$, we cut $S^1$ in four
intervals $I_{2,1},J_{2,1},I_{2,2},J_{2,2}$ of length $1/4$, each
one on the right of the previous one: $X_2$ is then defined to be the
map which linearly shifts these intervals over by 2, thus sending the
$I$'s onto the $I$'s and the $J$'s onto the $J$'s.  The map $X_2$ is thus the
rotation map by $\pi$.  Assume now we have built $X_n$ and we want to
build $X_{n+1}$.  Take the $2(n!)$ intervals of the partition
associated to $X_n$ and divide each of the intervals $I_{n,i}$
according to the proportions given by the partition $P_{n+1}$, cutting each $I_{n,i}$ into $2n+1=2(n+1)-1$ intervals.  Leave all of the $J_{n,i}$'s undivided.  We have
partitioned $S^1$ into
\[
n! + (2n+1) n! = 2[(n+1)!]
\]
intervals with dyadic endpoints. Starting with $J_{n+1,1}:=J_{n,1}$ we
relabel all the intervals of the new partition by $I$'s and $J$'s,
alternating them. The new piecewise linear map $X_{n+1}\in T$ is then
defined by shifting all the intervals by 2 (see figure
\ref{fig:Q-mod-Z-3-x} to see the construction of the maps $X_2$ and
$X_3$).
\begin{figure}[0.5\textwidth]
 \begin{center}
  \includegraphics[height=4cm]{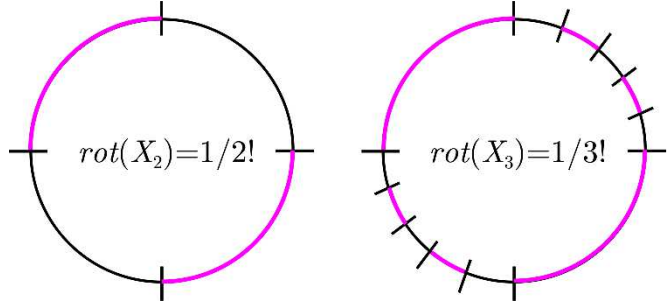}
 \end{center}
 \caption{Building the map $X_3$ from $X_2$.}
 \label{fig:Q-mod-Z-3-x}
\end{figure}
We need to verify that $(X_{n+1})^{n+1}=X_n$. We observe that
$Y_n:=(X_{n+1})^{n+1} \in T$ shifts every interval linearly by
$2n+2$. By construction $Y_n$ sends $J_{n,i}$ linearly onto
$J_{n,i+1}$, while it sends $I_{n,i}$ piecewise-linearly onto
$I_{n,i+1}$. All the possible breakpoints of $Y_n$ on the interval
$I_{n,i}$ occur at the points of the partition $P_{n+1}$, but it is
a straightforward computation to verify that the left and right
slope coincide at these points, thus showing that $Y_n$ sends
$I_{n,i}$ linearly onto $I_{n,i+1}$.

\medskip
\noindent \emph{Defining the embedding $\varphi$.}
To build the embedding $\varphi: \mathbb{Q}/\mathbb{Z} \to T$ we
define $\varphi(x_n):=X_n$ and extend it to a group homomorphism by recalling that $\mathbb{Q}/\mathbb{Z}=\langle x_n\rangle$.

The map $\varphi$ is easily seen to be injective.  If
$\varphi(x)=\mathrm{id}_{S^1}$ and $x=x_{i_1}^{m_{i_1}} \ldots
x_{i_\ell}^{m_{i_\ell}}$,
 then
\[
\mathrm{id}_{S^1} = X_{i_1}^{m_{i_1}} \ldots X_{i_\ell}^{m_{i_\ell}}.
\]
Since $(X_{r+1})^{r+1}=X_r$ for any integer $r$, we can rewrite the product $X_{i_1}^{m_{i_1}} \ldots X_{i_\ell}^{m_{i_\ell}}$
as $(X_n)^m$ for some suitable integers $n,m$. Since $\mathrm{id}_{S^1}=\varphi(x)=(X_n)^m$, we get that $m$ is a multiple of $n!$
and we can rewrite $x$ as $mx_n=(n!)x_n=0$.

\medskip
\noindent \emph{For every $x,y \in \mathbb{Q}/\mathbb{Z}, x \ne y$ the intervals
$\varphi(x)(J_{2,1})$ and $\varphi(y)(J_{2,1})$ are disjoint.}
If we define $V=\varphi(y)(J_{2,1})$, then the two intervals can be rewritten as $\varphi(x y^{-1})(V)$ and $V$.
Since $\varphi$ is an embedding and $xy^{-1} \ne 1$, these intervals must be distinct.
\end{proof}

As an immediate consequence of the previous proposition, we get the following two results from the introduction.

\vspace{.1 in}
\noindent{\bf Theorem 1.6} {\it For every $K \le \mathbb{Q} / \mathbb{Z}$ there is an embedding $F \wr_r K \hookrightarrow T$,
where $F$ and $T$ are the respective R. Thompson's groups.}

\begin{proof} We prove it for the full group $K = \mathbb{Q} / \mathbb{Z}$. We apply the previous
Theorem to build an embedding $\varphi:\mathbb{Q}/\mathbb{Z} \hookrightarrow T$. Moreover,
by construction, the image $\varphi(\mathbb{Q}/\mathbb{Z})$ acts as permutations on the intervals $\{J_{n,i}\}_{n,i \in \mathbb{N}}$.
Hence, we recover that
\[
\PL_2(J_{2,1}) \wr \mathbb{Q}/\mathbb{Z} \hookrightarrow T.
\]
where here $\PL_2(J_{2,1})$ is the subgroup of $T$ which consists of elements which are the identity off of $J_{2,1}$, that is, a group isomorphic with $F$.
\end{proof}

\vspace{.1cm}
\noindent{\bf Theorem 1.7} {\it For every $K \le \mathbb{Q} / \mathbb{Z}$
there is an embedding $\PL_+(I) \wr_r K \hookrightarrow \PL_+(S^1)$.}

\begin{proof} The proof of this result is similar to the one of Theorem \ref{thm:embedding-thompson},
except that here we do not require the endpoints of the interval $I$ to be dyadic.
\end{proof}

\remark{We remark that none of the proofs of the embedding results require the groups to have no
non-abelian free subgroups, although we notice that this condition
is automatically guaranteed in Theorems \ref{thm:embedding-thompson} and \ref{thm:embedding-PL}
because of the Brin-Squier Theorem (Theorem 3.1 in \cite{brin1}). However, in
Theorem \ref{thm:first-embedding} we may have non-abelian free subgroups inside
$H_0 \le \homeooI$ and still build the embedding.}

\appendix
\section{\label{sec:appendix}A counterexample to a construction of Solodov}

\begin{figure}[0.5\textwidth]
 \begin{center}
  \includegraphics[height=4.5cm]{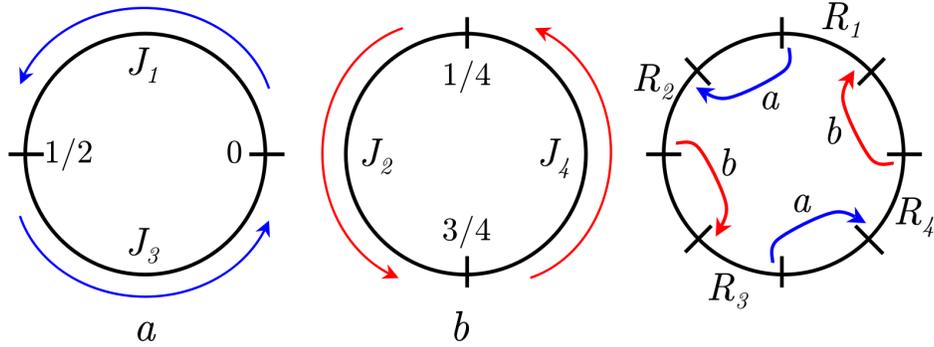}
 \end{center}
\caption{Construction and behavior of the maps $a$ and $b$.}
\label{fig:solodov-counter-example}
\end{figure}

Write the unit circle $S^1$ as $\mathbb{R}/\mathbb{Z}$ and define the intervals
\[
J_1= \left[0,\frac{1}{2}\right], \qquad
J_2=\left[\frac{1}{4},\frac{3}{4}\right], \qquad
J_3=\left[\frac{1}{2},0\right], \qquad
J_4=\left[\frac{3}{4},\frac{1}{4}\right]
\]
and the intervals
\[
R_1=\left[\frac{1}{8},\frac{1}{4}\right], \qquad
R_2=\left[\frac{3}{8},\frac{1}{2}\right], \qquad
R_3=\left[\frac{5}{8},\frac{3}{4}\right], \qquad
R_4=\left[\frac{7}{8},0\right].
\]
which are written making use of the local ordering of $S^1$. Finally, let 
$A_1,A_2,A_3,A_4$ be the left endpoints of $R_1,R_2,R_3,R_4$, that is
$A_1=\frac{1}{8}, A_2 = \frac{3}{8}, A_3 = \frac{5}{8}, A_4 = \frac{7}{8}$
(see figure \ref{fig:solodov-counter-example}).
We notice that
$A_1 \in J_4 \cap J_1, A_2 \in J_1 \cap J_2, A_3 \in J_2 \cap J_3, A_4 \in J_3 \cap J_4$.

We want to build maps $a,b$ in Thompson's group $T$ which act as in the figure.
Let $f:\left[0,\frac{1}{2}\right] \to \left[0,\frac{1}{2}\right]$ be
following piecewise linear map:
\[
f(t):=
\begin{cases}
4t & t \in \left[0,\frac{3}{32}\right] \\
t+\frac{9}{32} & t \in \left[\frac{3}{32}, \frac{1}{8}\right] \\
\frac{1}{4}t+\frac{3}{8} & t \in \left[\frac{1}{8},\frac{1}{2}\right].
\end{cases}
\]
It is immediate to verify that $f(R_1)=
\left[\frac{13}{32},\frac{7}{16} \right] \subseteq \ovc{R_2}$ 
(where $\ovc{R_2}$ denotes the interior of $R_2$).
Now we define $a$ to be equal to $f$
on $\left[0,\frac{1}{2} \right]$ and to ``act like $f$'' on $\left[\frac{1}{2},1 \right]$, that is
we conjugate $f$ by a rotation by half a circle
\[
a(t):=\begin{cases}
f(t) & t \in \left[0,\frac{1}{2} \right] \\
\rho_{\frac{1}{2}}f\rho_{\frac{1}{2}}^{-1}(t) & t \in \left[\frac{1}{2},1 \right],
\end{cases}
\]
where $\rho_{\frac{1}{2}}(t)=t+\frac{1}{2} \pmod{1}$ is the required rotation map.
The homeomorphism $a$ is in the standard copy of Thompson's group $F$, within Thompson's group $T$. To build $b$, we conjugate $a$ by a rotation by
a quarter of the circle. More precisely, we define
\[
b(t):=\rho_{\frac{1}{4}}a\rho_{\frac{1}{4}}^{-1}(t) \pmod{1}
\qquad t \in [0,1].
\]
where $\rho_{\frac{1}{4}}(t)=t+\frac{1}{4} \pmod{1}$ is the required rotation map.
The homeomorphism $b$ is in Thompson's group $T$.

We now follow Solodov's construction from the proof of Lemma 2.4 in \cite{solodov84} to find an element of $\langle a,b\rangle$ which has rotation number non-zero.  Of course, for this choice of $a$ and $b$, such elements are easy to find, but we are testing here the actual construction employed by Solodov.

By construction, we observe that $a$ and $b$ are in Thompson's group $T$ 
with $\Supp(a)=\ovc{J_1} \cup \ovc{J_3},
\Supp(b)=\ovc{J_2} \cup \ovc{J_4}$.
We also note that $a$ and $b$ satisfy the following inclusions (see figure
\ref{fig:solodov-counter-example})
\[
a(R_1) \subseteq \ovc{R_2}, \qquad  
b(R_2) \subseteq \ovc{R_3}, \qquad
a(R_3) \subseteq \ovc{R_4}, \qquad
b(R_4) \subseteq \ovc{R_1}.
\]
It is easy to see that our choice of $a,b$ yields that $\Supp \langle a,b \rangle=S^1$
and
$a(A_1)>A_2, b(A_2)>A_3, a(A_3)>A_4, b(A_4) > A_1$.

\medskip
Moreover, we observe that the element $baba$ (composing right-to-left) sends $A_1$ 
`Around the circle and past itself', which is effectively the condition created by Solodov's construction, 
and which is used to verify that such a constructed element would have non-zero rotation number.
However, we observe that $baba$ also sends the closed
interval $R_1$ inside the open interval $\ovc{R_1}$, therefore $baba$
fixes some point in $R_1$ by the contraction mapping lemma, and so $baba$
has rotation number zero. It can be verified that the attracting fixed
point for $baba$ detected by the previous argument is 
$\frac{1}{6} \in R_1$.

\bibliographystyle{plain}

 \small{ \noindent
  \textsc{
 \rule{0mm}{6mm} \\
 Collin Bleak \\
 School of Mathematics and Statistics, \\
 University of St Andrews, \\
 Mathematical Institute, \\
 North Haugh, \\
 St Andrews, Fife KY16 9SS, \\
 Scotland}
 \\ \emph{E-mail address:}
\texttt{collin@mcs.st-and.ac.uk}

 \small{ \noindent 
 \textsc{
 \rule{0mm}{6mm} \\
 Martin Kassabov \\
 Department    of Mathematics, Cornell University, \\
 590 Malott Hall, Ithaca, NY 14853,\\
 U.S.A.}
 \\ \emph{E-mail address:}
  \texttt{kassabov@math.cornell.edu}}

 \small{ \noindent 
 \textsc{
    \rule{0mm}{6mm} \\
 Francesco Matucci \\
Department of Mathematics, University of Virginia, \\
 325 Kerchof Hall, Charlottesville, VA 22904, \\
 U.S.A.}
 \\ \emph{E-mail address:} \texttt{fm6w@virginia.edu}}
\end{document}